\documentclass[10pt]{amsart}

\usepackage[english]{babel}
\usepackage[T1]{fontenc}
\usepackage{amsmath,amssymb,stmaryrd}
\usepackage[all]{xy}
\usepackage{graphicx}
\usepackage{fancyhdr}

\newtheorem{definition}{Definition }[section]

\newtheorem{lemma}[definition]
{Lemma}
\newtheorem{theorem}[definition]
{Theorem }
\newtheorem{ex}[definition]
{Example }

\newtheorem{corollary}[definition]
{Corollary }
\newtheorem{prop}[definition]{Proposition}

\newtheorem{rmk}[definition]{Remark}

\renewcommand{\ker}{\text{Ker}}

\newcommand{\lgw}{\longrightarrow}
\newcommand{\lgm}{\longmapsto}
\newcommand{\s}{\sigma}
\newcommand{\ovl}{\overline}
\newcommand{\Frac}{\text{Frac}}
\renewcommand{\deg}{\text{deg}\,}
\newcommand{\ord}{\text{ord}}
\newcommand{\Ker}{\text{Ker}}
\renewcommand{\Im}{\text{Im}}
\newcommand{\wdh}{\widehat}
\newcommand{\ini}{\text{in}}
\newcommand{\wdt}{\widetilde}
\renewcommand{\l}{\lambda}

\newcommand{\m}{\mathfrak{m}}

\newcommand{\Gr}{\text{Gr}}
\renewcommand{\k}{\Bbbk}
\renewcommand{\dim}{\text{dim}}
\newcommand{\trdeg}{\text{tr.deg}}
\newcommand{\cha}{\text{char\,}}
\newcommand{\R}{\mathbb{R}}
\newcommand{\K}{\mathbb{K}}
\newcommand{\grk}{\text{grk\,}}
\newcommand{\N}{\mathbb{N}}

\newcommand{\C}{\mathbb{C}}

\renewcommand{\t}{\tau}
\renewcommand{\a}{\alpha}
\renewcommand{\b}{\beta}
\newcommand{\g}{\gamma}

\renewcommand{\phi}{\varphi}

\newcommand{\e}{\varepsilon}
\begin{document}
\title[Homomorphisms of local algebras]{Homomorphisms of local algebras in positive characteristic\footnote{The author's research was partially supported by NSERC grant 0GP0009070}}

\author{Guillaume Rond}
\address{IML, Campus de Luminy, Case 907, 
13288 Marseille Cedex 9, France}
\email{rond@iml.univ-mrs.fr}
\maketitle

\begin{center}\textit{Dedicated to Professor Shuzo Izumi}\end{center}
\section{Introduction}
The aim of this paper is to study some properties  of regularity of homomorphisms of local  $\k$-algebras, in particular when $\k$ is a field of positive characteristic. In characteristic zero, the geometric rank of a homomorphism of local algebras $\phi : A\lgw B$ (denoted by $\grk(\phi)$) is a nice invariant that gives information about the structure of this homomorphism. In particular, a result due to P. M. Eakin and G. A. Harris \cite{E-H} asserts that a homomorphism between rings of formal power series (or convergent power series) over a field of characteristic zero can be monomialized, and after monomialization the geometric rank is equal to the dimension of the image of the monomial homomorphism. Homomorphisms with maximal geometric rank have nice properties that we can summarize in the following theorem:

\begin{theorem}\label{rappels}\cite{Ga2}\cite{E-H}\cite{Mi}\cite{B-Z}\cite{Iz3}
Let $\phi\ :\ A\lgw B$ be a homomorphism of analytic $\C$-algebras where $B$ is an integral domain. Then the following properties are equivalent:
\begin{itemize}
\item[i)] $\grk(\phi)=\dim(\wdh{A}/\Ker(\wdh{\phi}))$.
\item[ii)] $\grk(\phi)=\dim(\wdh{A}/\Ker(\wdh{\phi}))=\dim(A/\Ker(\phi))$.
\item[iii)] $\exists a\geq 1,\ b\geq 0$ such that $\wdh{\phi}^{-1}(\m_B^{an+b})\subset Ker(\wdh{\phi})+\m_A^n$ $\forall n\in\N$.
\item[iv)] $\wdh{\phi}(\wdh{A})\cap B=\phi(A)$.\end{itemize}
\end{theorem}
Moreover S. Izumi proved $ii)\Longleftrightarrow  iii)$ for any homomorphism of local rings of equicharacteristic zero \cite{Iz5}. \\

In characteristic zero, the geometric rank of $\phi :A\lgw B$ is equal to the rank of the $B$-module generated by $\Omega^1_{\k}(A)/\cap\m_A^n$ in $\Omega^1_{\k}(B)/\cap\m_B^n$. Unfortunately this definition  does not extend  well to positive characteristic for some obvious reasons (for instance look at the $\k$-homomorphism $\phi\ :\ \k[[x]]\lgw \k[[x]]$ defined by $\phi(x)=x^p$ where $\cha(\k)=p$: this homomorphism obviously satisfies a linear Chevalley estimate).
 
In this paper we  extend the definition of geometric rank in positive characteristic, using the transcendence degree of  the $\m_B$-adic valuation restricted to $A$ (cf. Section \ref{section-geo-rank}). This last  definition was first given by M. Spivakovsky in \cite{Sp}. In Section \ref{section-mono} we  prove a result (cf. Theorem \ref{structuretheorem}) about the structure of $\k$-homomorphisms between rings of power series over an infinite field of positive characteristic (similar to the result of P. M. Eakin and G. A. Harris \cite{E-H} valid in characteristic zero). This result involves our definition of geometric rank and shows that it is the right analogue of the geometric rank defined usually in characteristic zero. This result is very close to a monomialization result of dominant homomorphisms in positive characteristic. Moreover the proof of it is algorithmic and allows us to compute the geometric rank. \\
In Section \ref{section-Chevalley} we can deduce our first main result which is the positive characteristic analogue of the main result of \cite{Iz3}, i.e. $ii)\Longleftrightarrow iii)$ of Theorem \ref{rappels} (Linear Chevalley's Lemma):\\

\textbf{Theorem \ref{maintheorem}} 
\textit{ Let $\phi : A\lgw B$ be a homomorphism of local $\k$-algebras where $\k$ is a field of positive characteristic. Assume that $\wdh{A}$ is an integral domain and $B$ is regular. Then the following conditions are equivalent:
\begin{itemize}
\item[i)]$\grk(\phi)=\dim (A)$.
\item[ii)] There exist $a,b\in\R$ such that $a\nu_{\m_A}(f)+b\geq \nu_{\m_B}(\phi(f))$ for any $f\in A$.
\end{itemize}}
Homomorphisms satisfying these equivalent conditions are called \textit{regular} homomorphisms.

We would like to mention the work of R. H\"ubl \cite{Hu}  who gave sufficient conditions for general homomorphisms of local rings to satisfy condition $ii)$. He uses  a deep result of S. Izumi and D. Rees about the so-called Rees valuations. Unfortunately these conditions are  difficult to check in practice and we do not know if they are necessary conditions.\\
In characteristic zero, the result of Izumi is important in subanalytic geometry,  since Bierstone and Milman showed the paramount importance of the dependance of linearity of the Chevalley function on parameters for the composite function property (cf. \cite{B-M1}, \cite{B-M2} or see \cite{B-M3} for a general survey about the importance of the Gabrielov's Theorem and the Izumi's Theorem in subanalytic geometry).\\

The end of the paper is devoted to show how we can use the monomialization Theorem (Theorem \ref{structuretheorem}) in any charateristic in order to obtain new  results or generalizations of known results about regular homomorphisms of local $\k$-algebras in any characteristic. For example in the second part of Section \ref{section-Chevalley} we give an interpretation in terms of diophantine approximation of the fact that the Chevalley function of a homomorphism that is not regular is not bounded by an affine function.

In Section \ref{section-Henselian} we study  homomorphisms of Henselian $\k$-algebras, where $\k$ is a  field of any characteristic (for definitions, see Section \ref{section-Henselian}), which are generalizations of homomorphisms of convergent or formal power series rings, and we  give some cases where such a homomorphism $\phi : A\lgw B$ satisfies $\wdh{\phi}^{-1}(B)=A$ or $\Ker(\wdh{\phi})=\Ker(\phi)\wdh{A}$. For this we first state  a preparation theorem for Henselian $\k$-algebras (cf. Proposition \ref{Weierstrass}). Indeed the Weierstrass Preparation Theorem is essentially the only tool that we need for this study. Then we give a Henselian version of  Theorem \ref{structuretheorem} in any characteristic (cf. Theorem \ref{structuretheorem'}).  We deduce from it a weak version of a theorem of A. M. Gabrielov \cite{Ga2} (the analogue of $ii)\Longrightarrow iv)$ of Theorem \ref{rappels}) for good Henselian  $\k$-algebras in any characteristic (cf. Definition \ref{good_hens} for the definition of a good Henselian $\k$-algebra). 
This is our second main result:\\

\textbf{Theorem \ref{grk-strongly}.}
\textit{Let $\k$ be a field of any characteristic. Let $\phi : A\lgw B$ a homomorphism of good Henselian $\k$-algebras, where $A$ and $B$ are regular. If $\grk(\phi)=\dim (A)$ then $\phi$ is strongly injective (i.e. $\wdh{\phi}^{-1}(B)=A$).}\\

In Section \ref{section-ex} we study two particular cases of homomorphisms of local rings that are regular. First, using the algorithmic proof of Theorem \ref{structuretheorem'},  we prove that $\phi$ is injective if and only $\grk(\phi)=2$, when $A$ is a two dimensional integral $\k$-algebra with respect to a W-system (cf. Theorem \ref{dim2}). We deduce from this a generalization of a theorem due to S. S. Abhyankar and M. van der Put \cite{A-vdP} (who studied the case when $A$ is an analytic regular $\k$-algebra in any characteristic):\\

\textbf{Theorem \ref{A-vdP}.}
\textit{Let $\k$ be a field of any characteristic. Let $\phi : A\lgw B$ be a homomorphism of good Henselian  $\k$-algebras where $A$ is regular and $\dim(A)=2$.  If $\phi$ is injective then it is strongly injective.}\\

The second particular case is the case of homomorphisms of analytic algebras defined by algebraic power series over any valued field. This case has been previously studied for $\k=\C$ in \cite{To1}, \cite{Be} and \cite{Mi} using transcendental methods. We show here how to prove that such homomorphisms are regular using the monomialization theorem (cf. Corollary \ref{To}).\\
\\Finally, there are still remaining open problems. One of them is to know if the Gabrielov's Theorem (that asserts the following: if $\phi : A\lgw B$, a homomorphism of analytic $\C$-algebras, satisfies $\grk(\phi)=\dim(\wdh{A}/\Ker(\wdh{\phi}))$ then $\Ker(\wdh{\phi})=\ker(\phi)\wdh{A}$) extends to positive characteristic for analytic $\k$-algebras (and even for good Henselian $\k$-algebras in any characteristic). The proof of A. Gabrielov is quite difficult and the attempts to give a simpler proof (even over the field of complex number numbers $\C$) have not been successful. An other one is to extend these results in mixed characteristic. \\
\\
The author would like to thank Edward Bierstone for his comments about this work and his financial support at the University of Toronto. I would like also to thank  Mark Spivakovsky who pointed out a few mistakes in a preliminary version and helped me to improve this paper with  his mathematical and English remarks. Finally I would like to thank the referee for his useful comments and suggestions.

\subsection{Terminology }In this paper, rings are always assumed to be commutative Noetherian rings with unity. In any case $\k$ denotes a field. A \textit{local $\k$-algebra} will be a local ring $A$, with maximal ideal $\m_A$, along with an injective homomorphism $\k\lgw A$ such that the induced homomorphism $\k\lgw A/\m_A$ is a finite field extension. A \textit{homomorphism of local rings} $\phi\ :\ A\lgw B$ means a ring homomorphism such that $\phi(\m_A)\subset \m_B$ and the induced homomorphism  $A/\m_A\lgw B/\m_B$ is a finite extension of fields. The $\m_A$-adic order $\nu_{\m_A}$ is defined by $\nu_{\m_A}(f):=\max\{n\in\N\ /\ \ f\in\m_A^n\}$ for any $f\in A$.  For any $f\in A$, where $A$ is a local ring, $\ini(f)$ will denote the image of $f$ in $\Gr_{\m_A}A$.

\section{The geometric rank}\label{section-geo-rank}
Let $\phi : A\lgw B$ be a homomorphism of local $\k$-algebras and let us assume that $A$ is an integral domain and $B$ is regular.
Consider the valuation $\nu=\nu_B\circ\phi$ defined on $\Frac(A/\Ker(\phi))$, the quotient field of the domain $A/\Ker(\phi)$. We denote by $A_{\nu}$ the valuation ring associated to $\nu$ and by $\m_{\nu}$ its maximal ideal. We denote by $\trdeg_{\k}\nu$ the transcendence degree of the field extension $\k\lgw \frac{A_{\nu}}{\m_{\nu}}$.\\
The Abhyankar's Inequality says in our context that 
$$\trdeg_{\k}\nu+1\leq \dim (A/\ker(\phi))\ \ \  (\leq \dim(A)).$$

\begin{definition}\cite{Sp}
If $\Ker(\phi)\neq \m_A$, the integer $\trdeg_{\k}\nu+1$ is called the  \emph{geometric rank} of $\phi$ and denoted $\grk(\phi)$. If $\Ker(\phi)=\m_A$, then $\grk(\phi):=0$.
\end{definition}

\begin{lemma}\label{lemma3}
Let $\phi $, $A$ and $B$ as above. Assume moreover that $\wdh{A}$ is an integral domain. Then $\grk(\phi)=\grk(\wdh{\phi})$.
\end{lemma}

\begin{proof}
We denote by $A_{\nu}$ (resp. $\wdh{A}_{\wdh{\nu}}$) the valuation ring associated to $\nu=\nu_B\circ\phi$ (resp. to $\wdh{\nu}=\nu_B\circ\wdh{\phi}$) and $\m_{\nu}$ (resp. $\wdh{\m}_{\wdh{\nu}}$) its maximal ideal. We have $\wdh{\m}_{\wdh{\nu}}\cap A_{\nu}=\m_{\nu}$ thus the quotient homomorphism $\frac{A_{\nu}}{\m_{\nu}}\lgw \frac{\wdh{A}_{\wdh{\nu}}}{\wdh{\m}_{\wdh{\nu}}}$ is injective, hence $\grk(\wdh{\phi})\geq\grk(\phi)$.\\
On the other hand, if the images of $f_1,...,\,f_q\in \wdh{A}_{\wdh{\nu}}$ in  the field $\k_{\wdh{\nu}}=\frac{\wdh{A}_{\wdh{\nu}}}{\wdh{\m}_{\wdh{\nu}}}$ are algebraically independent over $\k$, then we can consider elements $f'_1,...,\,f'_q\in A_{\nu}$ such that $f'_i-f_i\in\wdh{\m}_{\wdh{\nu}}$. Thus the images of $f'_1,...,\,f'_q$ in $\k_{\wdh{\nu}}$ are algebraically independent over $\k$ because their images coincide with the images of $f_1,...,\,f_q$. Hence $\grk(\phi)=\grk(\wdh{\phi})$.
\end{proof}

\begin{lemma}\label{lemma4}
Let $\phi : A\lgw B$ be a homomorphism of local $\k$-algebras where $\wdh{A}$ is an integral domain and $B$ is regular. If $\grk(\phi)=\dim(A)$ then $\phi$ and $\wdh{\phi}$ are both injective. 
\end{lemma}

\begin{proof}
Since $\dim(A)=\dim(\wdh{A})$ and $\Ker(\phi)$ and $\Ker(\wdh{\phi})$ are prime ideals, the assertion  follows from the Abhyankar's Inequality.
\end{proof}

\begin{lemma}\label{lemma5}
Let $\phi : A\lgw B$ and $\sigma\ :\ A'\lgw A$ be homomorphisms of local  $\k$-algebras, where $A$ and $A'$ are integral domains and $B$ is regular. If $\s$ is finite and injective then $\grk(\phi\circ\s)=\grk(\phi)$.
\end{lemma}

\begin{proof}
We may assume that $\phi$ is injective by replacing the local $\k$-algebras $A$ and $A'$ by $A/\Ker(\phi)$ and $A'/\s^{-1}(\Ker(\phi))$ respectively.
We denote by $\nu$ and $\nu'$ the valuations induced by $\phi$ and $\phi\circ\s$ respectively. Let $f\in A_{\nu}$. Then there are $a_i\in A_{\nu'}'$ and $k\in\N$ such that 
$$a_0f^k+a_1f^{k-1}+\cdots+a_k=0$$
because $\Frac(A')\subset\Frac(A)$ is finite. We can assume that at least one of the $a_i$'s satisfies $\nu'(a_i)=0$ by dividing the last relation by an element $a_{i_0}$ satisfying $\nu'(a_{i_0})=\min_i\nu'(a_i)$. Then, if $\nu(f)=0$,  we see that the image of $f$ in $\k_{\nu}$ satisfies a non-trivial integral equation over $\k_{\nu'}$. Then
 the field extension $\k_{\nu'}\lgw \k_{\nu}$ is algebraic and $\grk(\phi)=\grk(\phi\circ\s)$.
\end{proof}

\begin{lemma}\label{lemma6}
Let $\phi  : \k[[x_1,...,\,x_n]]\lgw \k'[[y_1,...,\,y_m]]$ be a homomorphism of formal power series rings where $\k\lgw \k'$ is finite. Let $\phi_{\k'}$ denote the induced homomorphism $\k'[[x_1,...,\,x_n]]\lgw \k'[[y_1,...,\,y_m]]$. Then $\grk(\phi)=\grk(\phi_{\k'})$.
\end{lemma}

\begin{proof}
The homomorphism $\k[[x_1,...,\,x_n]]\lgw \k'[[x_1,...,\,x_n]]$ is finite and injective. Thus the result follows from Lemma \ref{lemma5}.
\end{proof}
Finally we give a combinatorial  characterization of the geometric rank. For any $f\in\k[[y_1,...,\,y_m]]$, we denote by $\ini(f)$ the form of lowest degree in the power series expansion of $f$. We define a total ordering $<$ on $\N^m$ in the following way: for any $\a,\,\b\in\N^m$, we say that $\a<\b$ if $(|\a|,\,\a_1,...,\,\a_m)<(|\b|,\,\b_1,...,\,b_m)$ for the left-lexicographic ordering, where $|\a|:=\a_1+\cdots+\a_m$. This ordering induces a monomial ordering on $\k[[y_1,...,\,y_m]]$. If $M=a_{\a}y^{\a}$ is a monomial, we define $\exp(M):=\a$.
For any $f\in\k[[y_1,...,\,y_m]]$, we define $\ini_{<}(f)$ to be the monomial of least order in the power series expansion of $f$ and $\exp(f):=\exp(\ini_{<}(f))$.

\begin{prop}\label{cone}
Let $\phi : A:=\k[[x_1,...,\,x_n]]\lgw B:=\k[[y_1,...,\,y_m]]$ be a homomorphism of formal power series rings. Let $C$ be the minimal cone of $\R^m$ containing  $\exp(\phi(f))$ for any $f\in\k[[x_1,...,\,x_n]]$. Then $\grk(\phi)=\dim (C)$.
\end{prop}

\begin{proof}
Let us  denote by $\ord$ the $(y_1,...,\,y_m)$-adic valuation on $B$ and $\nu$ the valuation on $A$ equal to $\ord\circ\phi$.\\
Let $\ovl{f}\in\k_{\nu}$ being the image of $f\in A_{\nu}$. We may write $f=\frac{g}{h}$ where $g,\,h\in\k[[x_1,...,\,x_n]]$ and $\nu(g)=\nu(h)$. The homomorphism $\phi$ induces an injection $\k_{\nu}\lgw\k_{\ord}=\k\left(\frac{y_1}{y_m},...,\,\frac{y_{m-1}}{y_m}\right)$, and the image of $\ovl{f}$ under this injection is just $\frac{\ini(\phi(g))}{\ini(\phi(h))}$.\\
Let us denote $B':=\k[\ini(\phi(f))]_{f\in \k[[x_1,...,\,x_n]]}$ and $\K':=\Frac(B')$. First we will prove that $\grk(\phi)=\dim(B')$.\\
We have $\dim(B')=\dim(\text{Spec}(B'))$. But we can look at $B'$ as a graded ring because any homogeneous component of any element of $B'$ is in $B'$. If we consider $\text{Proj}(B')$, we see that $\dim(B')=\dim(\text{Proj}(B'))+1$. So we have to prove that $\trdeg_{\k}\k_{\nu}$ is equal to the maximal number of algebraically independent elements of $\K'$ of the form $g/h$ where $g$ and $h$ are homogeneous of same degree.\\
Let us choose $f'_1,...,\,f'_r\in \K'$ algebraically independent over $\k$, such that $f'_i=\frac{g'_i}{h'_i}$ for any $i$, where $g'_i,\,h'_i\in B'$ are homogeneous of same degree.  By definition there exist $g_i$ and $h_i\in\k[[x_1,...,\,x_n]]$ such that $\ini(\phi(g_i))=g'_i$ and $\ini(\phi(h_i))=h'_i$ for any $i$. Let us denote $f_i:=\frac{g_i}{h_i}$ for any $i$. Then $f_i\in A_{\nu}$ for all $i$ and their images in $\k_{\nu}$ are algebraically independent over $\k$. Then we see that $\grk(\phi)\geq\dim(B')$.\\
On the other hand, let  $\ovl{f}_1,...,\,\ovl{f}_s\in \k_{\nu}$ be algebraically independent over $\k$. Let $f_i\in A_{\nu}$ be a lifting of $\ovl{f}_i$ for $1\leq i\leq s$. For any $i$ we may write $f_i=\frac{g_i}{h_i}$ where $g_i,\,h_i\in\k[[x_1,...,\,x_n]]$ and $\nu(g_i)=\nu(h_i)$. Let $f'_i$ denote $\frac{\ini(\phi(g_i))}{\ini(\phi(h_i))}$ for $1\leq i\leq s$. Then $f'_1,...,\,f'_s\in\K'$ are algebraically independent over $\k$. Thus $\grk(\phi)=\dim(B')$.\\
Now let us denote by $B''$ the sub-$\k$-algebra of $\k[y_1,...,\,y_m]$ generated by the $\ini_{>}(f)$ where $f\in B'$.
Then, because $B'$ is a $\k$-subalgebra of $\k[y_1,...,\,y_m]$ generated by homogeneous polynomials, the Hilbert function of $B'$ is the same as the Hilbert function of $B''$ (for instance look at Proposition 6.6.1 of \cite{K-R}). It implies that $\dim(B')=\dim(B'')$. But $\dim(B'')$ is exactly the dimension of $C$. Thus we have proved the proposition.
\end{proof}

\begin{corollary}\label{field_tr}
Let $\phi : A:=\k[[x_1,...,x_n]]\lgw B:=\k[[y_1,...,y_m]]$ be a homomorphism of formal power series rings. Let $t$ be a variable over $\k$ and let $\K:=\k(t)$. Let $\phi_{\K}: A':=\K[[x_1,...,\,x_n]]\lgw B':=\K[[y_1,...,\,y_m]]$ be the homomorphism of formal power series rings induced by $\phi$. Then $\grk(\phi_{\K})=\grk(\phi)$.
\end{corollary}

\begin{proof}

According to the proof of Proposition \ref{cone}, $\grk(\phi)=\dim(\k[\ini(f),\,f\in A])$ and  $\grk(\phi_{\K})=\dim(\K[\ini(f),\,f\in A'])$. If $g=\ini(f)$ with $f\in A'$ then $g=\ini(\l_1f_1+\cdots+\l_sf_s)$ with $\l_i\in\K$ and $f_i\in A$ for $1\leq i\leq s$. We may assume that $\l_i\in\k[t]$ for $1\leq i\leq s$ by multiplying $g$ by a non-zero element of $\K$. We write $\l_i=\sum_{j=0}^r\l_{i,j}t^j$ with $\l_{i,j}\in\k$ for $1\leq i\leq s$  and $0\leq j\leq r$. Then we get $g=\sum_{j=0}^rg_jt^j$ with $g_j=\ini(\sum_{i=1}^s\l_{i,j}f_i)$ for $0\leq j\leq r$ because $t$ is transcendant over $\k$. It follows that $\K[\ini(f),\,f\in A']$ and $ \K\otimes_{\k}\k[\ini(f),\,f\in A]$ are $\k$-isomorphic, hence $\grk(\phi)=\grk(\phi_{\K})$.
\end{proof}

\begin{prop}\label{compo}
Let $\phi :A\lgw B$, $\s_1\ :\ A\lgw A$ and $\s_2\ :\ B\lgw B$ be homomorphisms of local $\k$-algebras. Let us assume that there exist $a_1$, $a_2>0$ such that $\m_A^{a_1n}\subset \s_1(\m_{A}^n)$ and $\m_B^{a_2n}\subset \s_2(\m_{B}^n)$ for any $n\in\N$. Then $\grk(\s_1)=\dim(A)$, $\grk(\s_2)=\dim(B)$ and $\grk(\phi)=\grk(\phi\circ\s_1)=\grk(\s_2\circ\phi)$.
\end{prop}

\begin{proof}
We will prove the result for $\s_2$, the proof for $\s_1$ being similar. Using the notation used in the proof of Proposition \ref{cone}, $\grk(\phi)=\dim(B')$ is the degree of the Samuel polynomial $P(n)$ equal to $\dim_{\k}\left(\frac{\phi(A)}{\phi(A)\cap\m_B^n}\right)$ for $n>>0$. In the same way $\grk(\s_2\circ\phi)$ is equal to the degree of the Samuel polynomial $Q(n)$ equal to $\dim_{\k}\left(\frac{\s_2(\phi(A))}{\s_2(\phi(A))\cap\m_B^{n}}\right)$ for $n>>0$. By assumption we have $\m_B^{a_2n}\subset\s_2(\m_B^n)\subset\m_B^n$ for any $n\in\N$. Thus we get the following $\k$-linear maps:
$$\xymatrix{\frac{\s_2(\phi(A))}{\s_2(\phi(A))\cap\m_B^{a_2n}}\ar@{->>}[r] & \frac{\s_2(\phi(A))}{\s_2(\phi(A))\cap\s_2(\m_B^n)}\ar@{->>}[r] &\frac{\phi(A)}{\phi(A)\cap\m_B^n}\ar@{->>}[r] & \frac{\s_2(\phi(A))}{\s_2(\phi(A))\cap\m_B^n}}$$
where the first and last arrows are obvious quotient homomorphisms (thus they are $\k$-linear) and where the second arrow is a surjective $\k$-linear map defined by choosing a lifting in $\frac{\phi(A)}{\phi(A)\cap\m_B^n}$ of any element of $\frac{\s_2(\phi(A))}{\s_2(\phi(A))\cap\s_2(\m_B^n)}$.\\
Hence we have $Q(a_2n)\geq P(n)\geq Q(n)$ for $n>>0$. Thus we see that $\deg(P)=\deg(Q)$, hence $\grk(\s_2\circ\phi)=\grk(\phi)$. We get $\grk(\s_2)=\dim(B)$ by choosing $\phi=\text{id}_B$.
\end{proof}

\section{algorithm for modifying a homomorphism of a given rank}\label{section-mono}
We give here  a positive characteristic version of a theorem proved by Eakin and Harris \cite{E-H} in characteristic zero. This result is about the structure of homomorphisms of rings of formal power series over an infinite field of positive characteristic.
First we give the following definition:\\
\begin{definition}\label{modifications}
Let $\phi : \k[[x_1,...,\,x_n]]\lgw\k[[y_1,...,\,y_m]]$ a homomorphism of formal power series rings.
An \emph{admissible transformation} of $\phi$ is a homomorphism $\ovl{\phi}$ related to $\phi$ in one of the following ways:
\begin{enumerate}
\item \emph{Modification by automorphisms}: There exist a $\k$-automorphism $\t$ of $\k[[x_1,...,\,x_n]]$ and a $\k$-automorphism $\s$ of  $\k[[y_1,...,\,y_m]]$ such that $\ovl{\phi}=\s\circ\phi\circ\t$.
\item \emph{Modification by blowing-up}: There is $k\in \{1,...,\,m-1\}$ such that $\ovl{\phi}=\psi\circ\phi$ where $\psi$ is defined by
$$\psi(y_i)=y_i\ \ \text{ for } i\leq k,$$
$$\psi(y_i)=y_ky_i\ \ \text{ for } i>k.$$
\item \emph{Modification by ramification}: There is $d\in \N^*$ such that $\phi=\ovl{\phi}\circ\psi_d$ where $\psi_d$ is defined by 
$$\psi_d(x_1)=x_1^d, \text{ and } \ \psi_d(x_i)=x_i\ \ \forall i\neq 1.$$
\item \emph{Modification by contraction}: There is $k\in\{1,...,\,n-1\}$ such that $\phi=\ovl{\phi}\circ\psi$ where $\psi$ is defined by
$$\psi(x_i)=x_i\ \ \text{ for } i\leq k,$$
$$\psi(x_i)=x_ix_k\ \ \text{ for } i>k.$$
\end{enumerate}
\end{definition}
\begin{rmk}
We define the local $\k$-homomorphism $q\ :\ \k[[x_1,...,\,x_n]]\lgw \k[[x_1,...,,x_n]]$ by $q(x_1)=x_1x_2$  and $q(x_i)=x_i$ for $i>1$. It is clear that the homomorphisms $\psi$ defined in (4) of Definition \ref{modifications} are compositions of $q$ with permutations of the $x_i$'s. Thus we may use $q$ instead of $\psi$ in modification (4) of Definition \ref{modifications}. The same remark remains true for modifications by blowing-up.
\end{rmk}

\begin{lemma}\label{lemmamodifications}
Let $\phi : A:=\k[[x_1,...,\,x_n]]\lgw B:=\k[[y_1,...,\,y_m]]$ a homomorphism of formal power series rings. Let us consider a modification $\ovl{\phi}$ of $\phi$. Then $\grk(\ovl{\phi})=\grk(\phi)$. Moreover if there exist $a$ and $b$ such that 
$$a\nu_{\m_A}(f)+b\geq \nu_{\m_B}(\ovl{\phi}(f))$$ for any $f\in A$, then there exist $a'$ and $b'$ such that
$$a'\nu_{\m_A}(f)+b'\geq \nu_{\m_B}(\phi(f))$$ for any $f\in A$.
\end{lemma}

\begin{proof}
The lemma is obvious for modifications of type (1).\\
The second statement is a consequence of the following inequalities:
$$\nu_B(f)\leq\nu_B(\psi(f))\leq 2\nu_B(f) \ \ \forall f\in B, \ \text{ for modifications of type (2)},$$
$$\nu_A(f)\leq\nu_A(\psi_d(f))\leq d\,\nu_A(f)\ \ \forall f\in A, \ \text{ for modifications of type (3)},$$
$$\nu_A(f)\leq\nu_A(\psi(f))\leq 2\nu_A(f)\ \ \forall f\in A, \ \text{ for modifications of type (4)}.$$
Finally Proposition \ref{compo} gives us $\grk(\ovl{\phi})=\grk(\phi)$ in any cases.
\end{proof}

Now we can state the key result of this article. The proof of this theorem is inspired by the proof of a similar result in characteristic zero proved by Eakin and Harris \cite{E-H}.
\begin{theorem}\label{structuretheorem}Let $\k$ be an infinite  field of characteristic $p>0$ and consider a homomorphism $\phi : A:=\k[[x_1,...,\,x_n]]\lgw B:=\k[[y_1,...,\,y_m]]$ of power series rings.  Then there exists a finite sequence of admissible homomorphisms of formal power series rings $(\phi_i)_{i=0}^k\ :\ \k[[x_1,...,\,x_n]]\lgw\k[[y_1,...,\,y_m]]$ such that $\phi_0=\phi$ and $\phi_k(x_i)=y_i^{p^{\a_i}}u_i$, for some units $u_i$, for $i\leq \grk(\phi)$, and $\phi_k(x_i)=0$ for $i>\grk(\phi)$. Moreover, for any $i$, $u_i=1$ if $\a_i=0$, and $\ini(u_i)=1$ and $u_i\notin B^p$ if $\a_i>0$.
\end{theorem}
\begin{proof}
If $\grk(\phi)=0$, then $\phi(f)=0$ for all $f\in\m_A$. So we have the result.\\
We assume now that $\grk(\phi)>0$.\\
We will proceed by induction on the $q$-tuple $\mu=(\mu_1,...,\,\mu_q)\in (\N\cup\{+\infty\})^q$, defined later,  ordered with the lexicographic order where $q\leq n$. At the beginning, $q=n$ and this $q$-tuple is $(+\infty,...,\,+\infty)$.\\

\textbf{Step 0:} If $\phi(x_1)=0$  we exchange $x_n$ and $x_1$. Then we define $q=n-1$ and $\mu:=(\mu_1,...,\,\mu_{n-1})=(+\infty,...,\,+\infty)$.\\

\textbf{Step 1:} If $\phi(x_1)\neq0$, then we denote $d:=\ord(\phi(x_1))\in\N^*$. We denote by $g_d(y_1,...,\,y_m)$ the initial term of $\phi(x_1)$.  Let $(a_{i,j})_{i,j=1,...,\,m}$ be a non-singular matrix with entries in $\k$ such that $g_d(a_{1,\,1},\,...,\,a_{m,1})\neq 0$ ($\k$ is infinite). We define an automorphism $\psi$ of $\k[[y_1,...,\,y_m]]$ by
$$\psi(y_j):=\sum_{k=1}^{m}a_{j,\,k}y_k,\ \ \forall j=1,...,\,m.$$
So we get 
$$\psi\circ\phi(x_1)=g_d(a_{1,\,1},\,...,\,a_{m,1})y_1^d+\{\text{terms of degree $d$ not divisible by $y_1^d$}\}$$
$$+\{\text{terms of degree $>d$}\}.$$
By composing $\psi\circ\phi$ on the right by the automorphism of $\k[[x_1,...,\,x_n]]$ consisting in dividing $x_1$  by $g_d(a_{1,\,1},\,...,\,a_{m,1})$, we may assume that 
$$\phi(x_1)=y_1^d+\{\text{terms of degree $d$ not divisible by $y_1^d$}\}$$
$$+\{\text{terms of degree $>d$}\}.$$
Now we define the  homomorphism $\psi$  by
$$\psi(y_1):=y_1$$
$$\psi(y_i):=y_1y_i,\ \text{for }i>1.$$
We have $\psi\circ\phi(x_1)=uy_1^d$, $u$ being a unit of $\k[[y_1,...,\,y_m]]$ with $\ini(u)=1$.\\

\textbf{Step 2:} If $d=ep^{\a}$ with $\gcd(e,\, p)=1$, then we see that
$\psi\circ\phi=\phi'\circ\t'$ where $\t'(x_1)=x_1^e$ and $\t'(x_i)=x_i$ for $i\neq 1$, and $\phi'(x_1)=u'y_1^{p^{\a}}$, $\phi'(x_i)=\psi\circ\phi(x_i)$ for $i\neq 1$ and $\ini(u')=1$. So we can replace $\phi$ by $\phi'$.\\
In particular, if $\gcd(d,\, p)=1$, then we may assume $\phi(x_1)=y_1$.\\
Then, if $u\in B^{p^\b}$, with $\b\leq \a$, then we have $\phi(x_1)=u'^{p^{\b}}(y_1^{p^{\a-\b}})^{p^{\b}}$. So we see that
$\psi\circ\phi=\phi'\circ\t'$ where $\t'(x_1)=x_1^{p^{\b}}$ and $\t'(x_i)=x_i$ for $i\neq 1$, and $\phi'(x_1)=u'y_1^{p^{\a-\b}}$, $\phi'(x_i)=\psi\circ\phi(x_i)$ for $i\neq 1$ and $\ini(u')=1$.\\
So we may assume that $\phi(x_1)=uy_1^{p^{\a_1}}$, $\ini(u)=1$ and $u\notin B^p$ if $\a_1\neq 0$. At this step, the $q$-tuple $(\mu_1,...,\,\mu_q)=(\a_1,\,+\infty,...,\,+\infty)$.\\

\textbf{Step 3:}
We assume that $\phi(x_i)=y_i^{p^{\a_i}}u_i$, for $i<j$, with $\ini(u_i)=1$ and $u_i\notin B^p$ if $\a_i\neq 0$, and $\phi(x_i)=0$ for $i>q$. Moreover we assume that $\a_1\leq\a_2\leq\cdots\leq\a_{j-1}$. We denote $\mu=(\a_1,...,\,\a_{j-1},\,+\infty,...,\,+\infty)\in(\N\cup\{+\infty\})^q$.\\
We assume that $\ini(\phi(x_j))$ contains a monomial of the form $cy_1^{k_1}...y_{j-1}^{k_{j-1}}$. If $p^{\a_i}$ divides $k_i$ for all $i\leq j-1$, then we replace $x_j$ by the element $x_j-cx_1^{k_1/p^{\a_1}}...x_{j-1}^{k_{j-1}/p^{\a_{j-1}}}$. We can go on and by induction, there are two cases. In the first case we can replace $x_j$  by an element of the form $x_j-\sum_{\underline{k}}c_{\underline{k}}x_1^{k_1/p^{\a_1}}...x_{j-1}^{k_{j-1}/p^{\a_{j-1}}}$, where the sum is finite, and then we can assume that $\ini(\phi(x_j))$ has no monomial of the form $cy_1^{k_1}...y_{j-1}^{k_{j-1}}$ where $p^{\a_i}$ divides $k_i$ for all $i$. In the second case we can replace $x_j$  by an element of the form $x_j-\sum_{\underline{k}}c_{\underline{k}}x_1^{k_1/p^{\a_1}}...x_{j-1}^{k_{j-1}/p^{\a_{j-1}}}$, where the sum is not necessarily finite, and then we have $\phi(x_j)=0$.\\
If $\phi(x_j)=0$, then we exchange $x_{q}$ and $x_j$. Then we replace $q$ by $q-1$ and $\mu:=(\mu_1,...,\,\mu_{q})=(\a_1,...,\,\a_{j-1},\,+\infty,...,\,+\infty)$.\\

\textbf{Step 4:}
We assume that $\phi(x_i)=y_i^{p^{\a_i}}u_i$, for $i<j$, for some units $u_i$ with $\ini(u_i)=1$ and $u_i\notin B^p$ if $\a_i\neq 0$, and $\phi(x_i)=0$ for $i>q$. Moreover we assume that $\a_1\leq\a_2\leq\cdots\leq\a_{j-1}$. As before we denote $\mu=(\a_1,...,\,\a_{j-1},\,+\infty,...,\,+\infty)\in(\N\cup\{+\infty\})^q$.\\
Let us consider $cy_1^{k_1}...y_m^{k_m}$ a monomial of $\ini(\phi(x_j))$. If one of $k_j,...,\,k_m$ is different from zero, then after a permutation of the elements $y_j,...,\,y_m$, we can assume that $k_j\neq 0$. According to Step 3 we can assume that $\ini(\phi(x_j))$ has no monomial of the form $cy_1^{k_1}...y_{j-1}^{k_{j-1}}$ where $p^{\a_i}$ divides $k_i$ for all $i$, and we assume that $\phi(x_j)\neq 0$.\\
Assume that for any non-zero monomial $M=cy_1^{k_1}...y_{j-1}^{k_{j-1}}$ of $\ini(\phi(x_j))$,  $p^{\a_i}$ divides $k_i$ for any $i<l$, but for at least  one such monomial $p^{\a_l}$ does not divide $k_l$. 
After a change of variables of the form $\s(y_i)=y_i$ for $i\leq l$ and $\s(y_i)=y_i+y_l^{\l_i}$ for $i>l$ and for some $\l_i\in\N$, we may assume that $\ini(\phi(x_j))$ contains a non-zero monomial of the form $cy_1^{k_1}...y_{l}^{k_l}$ where  $p^{\a_i}$ divides $k_i$ for any $i<l$ and $p^{\a_l}$  does not divide $k_l$. Then after a composition with a homomorphism of the form $\psi(y_i)=y_i$ for $i\leq l$ and $y_i=y_1y_i$ for $i>l$, we may assume that each monomial of $\ini(\phi(x_j))$ depends only on $y_1,...,\,y_l$. And by Step 3, we may assume that  for any monomial $cy_1^{k_1}...y_{l}^{k_l}$ of $\ini(\phi(x_j))$,  $k_i$ is divisible by $p^{\a_i}$ for any $i<l$,  and that $k_l$ is  not divisible by $p^{\a_l}$. Finally we can exchange $x_j$ and $x_l$ and we can apply the following lemma with $\a=\a_l$ :

\begin{lemma}\label{monomialization2}
Under the hypothesis of Theorem \ref{structuretheorem}, we assume that $\phi(x_i)=y_i^{p^{\a_i}}$ for all $i<l$ and that the monomials of $\ini(\phi(x_l))$ depend only on $y_1,...,\,y_l$. We assume moreover that for any monomial $cy_1^{k_1}...y_{l}^{k_l}$ of $\ini(\phi(x_l))$,  $k_i$ is divisible by $p^{\a_i}$ for any $i<l$,  and that $k_l$ is  not divisible by $p^{\a}$. Then there exists a finite sequence of modifications of $\phi$, such that $\ovl{\phi}$, the last homomorphism of the sequence, satisfies
$$\ovl{\phi}(x_i)=y_i^{p^{\a_i}}u_i\ \text{ for } i<l,$$
$$\ovl{\phi}(x_l)=y_l^{p^{\a'}}u_l$$
for some units $u_j$ and with $\a'< \a$.
\end{lemma}
\begin{proof}[Proof of Lemma \ref{monomialization2}]
We have $\phi(x_l)=M_1v_1+\cdots+M_rv_r$ for some units $v_i$ and some monomials $M_i$. We assume that this  expression is minimal: it means that none of these monomials divides another one. The convex hull in $\N^m$ of the set of elements $(w_1,...,\,w_m)$ such that $\ini(\phi(x_l))$ contains a non-zero monomial of the form $cy_1^{w_1}...y_m^{w_m}$ is a convex polyhedron $P$ of dimension strictly less than $l$ (because all such elements satisfy $w_{l+1}=\cdots=w_m=0$). Let $(w_1,...,\,w_l,\,0,...,\,0)$ be a vertex of this polyhedron. In particular $w_l$ is not divisible by $\a$. We may assume that $M_1$ corresponds to this vertex. We denote by $(w_{1,\,k},...,\,w_{m,\,k})$  the element of $\N^m$ that corresponds to $M_k$ for $k>1$. Because $M_1$ is a vertex of $P$ the cone defined by the following equations in the variables $e_i$:
$$\sum_{i=1}^l (w_{i,\,k}-w_i)e_i>0$$
for all $k$ such that the monomial $M_k$ depends only on $y_1,...,\,y_l$,
is a non-empty open set of $(\R_{\geq 0})^l$. Moreover for any $C>0$, using modifications by blowing-up on the variables $y_{l+1},...,\,y_m$, we may assume that the monomials $M_k$ depending on at least one of $y_{l+1},...,,\,y_m$ satisfy $w_{1,\,k}+\cdots+w_{l,\,k}>C$. Hence the cone defined by the equations:
\begin{equation}\label{cone1}\sum_{i=1}^l (w_{i,\,k}-w_i)e_i>0,\ \ k=2,...,\,r\end{equation}
is a non-empty open set of $(\R_{\geq 0})^l$.
Let  $(e_1,...,\,e_l)$ be $l$ linearly independent vectors of this cone with coefficients in $\N$: we write $e_i=(e_{i,\,1},...,\,e_{i,\,l})$ for each $i$. We can choose these vectors such that their determinant is not divisible by $p$ and such that $p$ does not divide $e_{l,\,l}$.\
Next we consider $\psi$ defined by:
$$\psi(y_i)=y_1^{e_{i,\,1}}...y_l^{e_{i,\,l}}\ \text{ for } 1\leq i\leq l$$
$$\psi(y_{i})=y_i\ \text{ for } i>l.$$
Hence, because the vectors $e_i$ satisfy (\ref{cone1}), $\psi\circ\phi(x_l)$ is of the form $\psi(M_1)u_l$ for some unit $u_l$. More precisely we have:
$$\psi\circ\phi(x_i)=y_1^{p^{\a_i}e_{i,\,1}}...y_l^{p^{\a_i}e_{i,\,l}}u_i\ \text{ for } i<l$$
$$\psi\circ\phi(x_l)=y_1^{\sum_{i=1}^lw_ie_{i,\,1}}...y_l^{\sum_{i=1}^lw_ie_{i,\,l}}u_l.$$
Because the vectors are linearly independent and because their determinant is not divisible by $p$, we can reduce to the following case by using modifications of type $(4)$:
$$\psi\circ\phi(x_i)=y_i^{p^{\a_i}}u_i\ \text{ for } i<l$$
$$\psi\circ\phi(x_l)=y_l^{w_le_{l,\,l}\det(e_{i,\,k})}u_l$$
for some units $u_i$. And because $e_{l,\,l}$ and $\det(e_{i,\,k})$ are not divisible by $p$, according to the Cramer's rule and using Step 2, we may assume that:
$$\phi(x_l)=y_l^{p^{\a'}}u'_l$$
where $\a'< \a$, because $w_l$ is not divisible by $p^{\a}$.
\end{proof}
Finally, using Step 2, we may assume that
$$\phi(x_i)=y_i^{p^{\a'_i}}u_i\ \text{ for } i\leq l,$$ 
for some units $u_i$,
where $(\a'_1,...,\,\a'_l)<_{lex}(\a_1,...,\,\a_l)$ (this can be achieved by permuting the $x_i's$ and the $y_i's$). Then, if we denote $\mu'=(\a'_1,...,\,\a'_l,\,+\infty,...,\,+\infty)$, we have $\mu>_{lex}\mu'$.\\

\textbf{Step 5:}
We assume that we have $\phi(x_i)=y_i^{p^{\a_i}}u_i$ for some units $u_i$ with $\ini(u_i)=1$, where $u_i\notin B^p$ whenever $\a_i>0$, for any $i<j$, and $\phi(x_i)=0$ for $i>q$. We assume that $\a_1\leq\a_2\leq\cdots\leq\a_{j-1}$ and we denote  $\mu=(\a_1,...,\,\a_{j-1},\,+\infty,...,\,+\infty)\in(\N\cup\{+\infty\})^q$.\\
According to Step 4, we may assume that none of the monomials of $\ini(\phi(x_j))$ depends only on $y_1,...,\,y_{j-1}$. After a change of variables in $y_j,...,\,y_m$ we may assume that one of the monomials of $\ini(\phi(x_j))$ depends only on $y_1,...,\,y_j$. By composing with the homomorphism $\psi$ defined by
$$\psi(y_i)=y_i,\ \text{ for } i\leq j$$
$$\psi(y_i)=y_iy_j,\ \text{ for } i>j$$
we can assume that $\ini(\phi(x_j))$ depends only on $y_1,...,\,y_j$, but any of its monomials depends on $y_j$. So we have $\ini(\phi(x_j))=y_j^{k}P_{d-k}(y_1,...,\,y_j)$ where $P_{d-k}$ is a homogeneous polynomial of degree $k$ not divisible by $y_j$.\\
Thus, we may use Lemma \ref{monomialization2} and assume that
$$\phi(x_i)=y_i^{\a_i}u_i\ \text{ for }i<j,$$
$$\phi(x_j)=y_j^{\a}u_j$$
for some units $u_i$ and some integer $\a$.
Finally, using Step 2, we may assume that $\phi(x_i)=y_i^{p^{\a'_i}}u_j$ for some units $u_i$ with $\ini(u_i)=1$ and $u_i\notin B^p$ if $\a'_i\neq 0$. Moreover we see that $(\a'_1,...,\,\a'_{j-1})\leq_{lex} (\a_1,...,\,\a_{j-1})$. Hence after  permutation of the variables we may assume that $\a'_1\leq\cdots\leq\a'_j$ and $\mu'=(\a'_1,...,\,\a'_j,\,+\infty,...,\,+\infty)<\mu$.\\

\textbf{Step 6:}
Eventually, we have $\phi(x_i)=y_i^{p^{\a_i}}u_i$, for $i\leq q$ and $\phi(x_i)=0$ for $i>q$, where $\a_1\leq \a_2\leq\cdots\leq \a_q$ and $u_i$ are units. In this case one checks that
$$\frac{A_{\nu}}{\m_{\nu}}=\k\left(\frac{x_q}{x_{q-1}^{p^{\a_q-\a_{q-1}}}},...,\,\frac{x_2}{x_1^{p^{\a_2-\a_1}}}\right),$$ and we have
$\grk(\phi)=q$. Because the geometric rank is invariant under modifications, we get the result.
\end{proof}

\begin{rmk}\label{char0}
If $\cha(\k)=0$, the proof of the result of Eakin and Harris is similar. Namely, at Step 2 we get $\phi(x_1)=y_1$ because any unit $u$ with $\ini(u)=1$ is a $d$-power for any $d\in\N^*$. Then we can skip Steps 3 and 4 because if $g(y_1,...,\,y_{j-1}):=\phi(x_j)(y_1,...,\,y_{j-1},\,0,...,\,0)$, then we replace $x_j$ by $x'_j=x_j-g(x_1,...,\,x_{j-1})$ and $\phi(x'_j)$ has no monomial of the form $cy_1^{k_1}...y_{j-1}^{k_{j-1}}$.
\end{rmk}

\section{Linear Chevalley's lemma}\label{section-Chevalley}
The aim of this section (and originally of the present paper) is to give an answer to a question that S. Izumi asked the author. This question is related to the following  result of C. Chevalley on complete local rings:
\begin{theorem}\cite{Ch}
Let $A$ be a complete local ring with maximal ideal $\m$. Let $(\mathfrak{a}_n)$ be a decreasing sequence of ideals of $A$ such that $\cap_n\mathfrak{a}_n=\{0\}$. Then there exists a function $\b : \N\lgw \N$ such that $\mathfrak{a}_{\b(n)}\subset \m^n$ for any positive integer $n$. 
\end{theorem}
In particular, if we consider an injective homomorphism of local rings $\phi : A\lgw B$ where $A$ is complete, then there exists a function $\b : \N\lgw \N$ such that $\phi^{-1}(\m_B^{\b(n)})\subset \m_A^n$ for any natural number $n$. This can be restated by saying that $\b(\nu_{\m_A}(f))\geq \nu_{\m_B}(\phi(f))$ for any $f\in A$. The least function $\b$ satisfying this inequality is called the \textit{Chevalley function} of $\phi$. If $\b$ is bounded from above by a linear function we say that $\phi$ has a linear Chevalley estimate.\\
S. Izumi (in \cite{Iz3} and \cite{Iz4}), in the case of equicharacteristic zero local rings, proved that $\phi$ has a linear Chevalley estimate if and only if $\grk(\phi)=\dim(A)$. The question asked by S. Izumi was the following: is it possible to extend this result for local $\k$-algebras with $\cha(\k)>0$? 
We can state such an analogue of the main result of \cite{Iz3} in positive characteristic :

\begin{theorem}\label{maintheorem}
Let $\phi : A\lgw B$ be a homomorphism of local $\k$-algebras where $\k$ is a field of positive characteristic. Assume that $\wdh{A}$ is an integral domain and $B$ is regular. Then the following conditions are equivalent:\\
\begin{itemize}
\item[i)]$\grk(\phi)=\dim (A)$.
\item[ii)] There is $a,b\in\R$ such that $a\nu_{\m_A}(f)+b\geq \nu_{\m_B}(\phi(f))$ for any $f\in A$.\\
\end{itemize}
\end{theorem}

\begin{definition}
A homomorphism of $\k$-algebras $\phi : A\lgw$ such that $\grk(\phi)=\dim(A/\Ker(\phi))$ is called a \textit{regular} homomorphism of $\k$-algebras.
\end{definition}\subsection{Proof of Theorem \ref{maintheorem}}

In order to give a proof of this theorem, we first state the following two lemmas:
\begin{lemma}\label{finite}
Let $\s : A\lgw B$ be a finite and injective homomorphism of complete local rings (we do not assume that the rings are local rings of equal characteristic). Then $\s$ satisfies property $ii)$ of Theorem \ref{maintheorem}.
\end{lemma}
\begin{proof}
By induction we only need to prove the lemma when $B$ is generated by a single element over $A$. We denote by $z$ this element which is integral over $A$. If $z\notin \m_B$, then $\m_B=\m_AB$, thus, for any $n\in\N^*$, $\m_B^n\cap A=\m_A^nB\cap A\subset \m_A^{n-C}$ for some $C\in\N$ not depending on $n$ (by the Artin-Rees Lemma).\\
Let us assume from now on that $z\in\m_B$. The $\k$-algebra $B/\m_AB$  is a finite $\k$-module generated by $1$,..., $z^{d-1}$ modulo $\m_A$ for some $d\in\N^*$. Let us assume that $1$,..., $z^{d-1}$ is a $\k$-basis of this $\k$-algebra. Thus, by Theorem 30.6 \cite{Na}, we see that $1$, $z$,..., $z^{d-1}$ generate $B$ as a $A$-module. It means that there exists an irreducible polynomial  $P(Z):=Z^d+a_{d-1}Z^{d-1}+\cdots+a_0\in A[Z]$ such that $B$ is isomorphic to $A[Z]/((P(Z))$. Moreover $a_i\in\m_A$, for $0\leq i\leq d-1$, because $1$, ..., $z^{d-1}$ is a $\k$-basis of $B/\m_AB$.
Let $\displaystyle\a:=\min_{0\leq i\leq d-1}\{\ord_A(a_i)\}$ (in particular $\a>0$). Then $z^d\in\m_A^{\a}B$. By induction, $z^{dn}\in\m_A^{\a n}B$ for any $n\in\N$. Hence, for any $n\in\N$:
$$\m_B^{(\a+d)n+1}=\left(\m_AB+(z)\right)^{(\a+d)n+1}\subset\m_A^{(\a+d)n+1}B+(z)\m_A^{(\a+d)n}B+\cdots+(z^{dn+1})\m_A^{\a n}B+(z^{dn})$$
$$\subset \m_A^{\a n}B+(z^{dn})\subset \m_A^{\a n }B.$$
Thus $\m_B^{(\a+d)n+1}\cap A\subset  \m_A^{\a n }B\cap A\subset \m_A^{\a n-C}$ for some $C\in\N$ not depending on $n$. This proves the lemma because $\a>0$.\\
\end{proof}

\begin{lemma}\label{redfinite}\cite{Iz1}
Let $\phi : A\lgw B$ and $\s : A'\lgw A$ be two homomorphisms of local rings where $\s$ is finite and injective and $\wdh{A}$ is an integral domain. Then $\phi$ satisfies $ii)$ if and only if $\phi\circ\s$ satisfies $ii)$.
\end{lemma}
\begin{proof}
Because $\s$ is finite and injective, if $\phi$ satisfies $ii)$ then $\phi\circ\s$ satisfies $ii)$ by Lemma \ref{finite}. 
In order to prove the "if"-part we follow the proof of Theorem 1.2 (1) $\Longrightarrow$ (2) of \cite{Iz1} using the fact that there exist two positive constants $c,\,d$ such that
$\nu_A(fg)\leq c(\nu_A(f)+\nu_A(g))+d,\ \forall f,\,g\in A$ (cf. \cite{Re})

\end{proof}

Now we can begin the proof of  Theorem \ref{maintheorem}. We will first reduce to the case where $A$ and $B$ are complete. We remark that for any $f\in A$ we have $\nu_{\m_A}(f)=\nu_{\m_{\wdh{A}}}(f)$ because $\wdh{A}$ is flat over $A$. So the order is invariant under completion. The Krull dimension and the geometric rank are also invariant under completion. Moreover, the inequality $ii)$ of Theorem \ref{maintheorem} is equivalent to a similar estimate for $\wdh{\phi} : \wdh{A}\lgw \wdh{B}$ following remark 4.4 of \cite{Iz1}.\\

From now on we assume that $A$ and $B$ are complete and $\phi=\wdh{\phi}$. Then we show that $i)$ and $ii)$ are always  true if $\dim(A)=0$: $\grk(\phi)=\dim(A)$ is trivially true because $\grk(\phi)\leq \dim(A)=0$. So $i)$ is true. In particular, using Lemma \ref{lemma4}, $\phi$ is injective. On the other hand $A$ is artinian and so the following descending chain of ideals stabilizes: $A\supset \phi^{-1}(\m_B)\supset...\supset\phi^{-1}(\m_B^n)\supset...$. So there exists $b$ such that $\nu_{\m_B}(\phi(f))\leq b$ for any $f\in A\backslash\Ker\,\phi=A\backslash\{0\}$, and $ii)$ is true.\\

From now on we assume that $A$ and $B$ are complete, $B$ is regular and $\dim(A)\geq 1$. In particular $B=\k'[[y_1,...,\,y_m]]$ where $\k\lgw\k'$ is finite (it follows from the definition of a homomorphism of local $\k$-algebras). \\

\textbf{Step 1:}
Assume that $\k=\k'$ and $\k$ is an infinite field.\\

(I) Implication $ii)\Longrightarrow i)$.\\
 We may reduce to the case $A$ is regular by using Lemma \ref{lemma5}, i.e. $A=\k[[x_1,...,\,x_n]]$. Moreover we have $B=\k[[y_1,...,\,y_m]]$. We need to prove that $\grk(\phi)=n$.\\
 Using Theorem \ref{structuretheorem}, we get a commutative diagram as follows:
 $$\xymatrix{\k[[x_1,...,\,x_n]]\ar[r]^{\phi} \ar[d]^{\s_1} & \k[[y_1,...,\,y_m]]\ar[d]^{\s_2}\\
 \k[[x_1,...,\,x_n]]\ar[r]^{\ovl{\phi}} & \k[[y_1,...,\,y_m]] }$$
 where $\s_1$ and $\s_2$ are compositions of homomorphisms defined in Definition \ref{modifications} and $\ovl{\phi}(x_i)=y_i^{p^{\a_i}}u_i$ for $1\leq i\leq \grk(\phi)$ and $\ovl{\phi}(x_i)=0$ if $i>\grk(\phi)$. Let us denote $r:=\grk(\phi)$. In particular we have $n\geq r$. By Lemma \ref{lemmamodifications}, we see that $\s_2\circ\phi$ satisfies property $ii)$. If $f\in\k[[x_1,...,\,x_n]]$ then $\s_2\circ\phi(f)=\ovl{\phi}\circ\s_1(f)$, thus the homogeneous component of minimal degree in the Taylor expansion of $\s_2\circ\phi(f)$ depends only on $y_1$,..., $y_r$. Thus, if we denote by $\pi\ :\ \k[[y_1,...,\,y_m]]\lgw \k[[y_1,...,\,y_r]]$ the canonical projection, we see that the order of $\s_2\circ\phi(f)$ is the same as the order of $\pi\circ\s_2\circ\phi(f)$ for any $f\in\k[[x_1,...,\,x_n]]$. Thus $\pi\circ\s_2\circ\phi$ satisfies property $ii)$. Moreover $\grk(\pi\circ\s_2\circ\phi)=\grk(\pi\circ\ovl{\phi}\circ\s_1)=r=\grk(\phi)$. Thus we may assume that $A=\k[[x_1,...,\,x_n]]$, $B=\k[[y_1,...,\,y_m]]$ and $\grk(\phi)=m\leq n$.\\
By assumption there exist $a$ and $b$ such that $\phi^{-1}(\m_B^{ak+b})\subset \m_A^k$ for any $k\in\N$. So, for any $k\in \N$, we may define surjective $\k$-linear maps
$$\phi(A)/\left(\m_B^{ak+b}\cap\phi(A)\right)\lgw A/\m_A^k$$
by choosing a lifting in $A/\m_A^k$ of any element of $\phi(A)/\left(\m_B^{ak+b}\cap\phi(A)\right)$.\\
Because $\phi(A)/\left(\m_B^{ak+b}\cap\phi(A)\right)$ is a $\k$-subspace of $B/\m_B^{ak+b}$, we have the following equalities and inequalities for any $k\in\N$:
$$(ak+b+m-1)!/\left((ak+b-1)!\,m!\right)=\dim_{\k}B/\m_B^{ak+b}$$
$$\geq \dim_{\k}\phi(A)/\left(\m_B^{ak+b}\cap\phi(A)\right)
\geq \dim_{k} A/\m_A^{k}=(n+k-1)!/\left((k-1)!\,n!\right)$$
Hence, by comparing the degree in $k$ of these two polynomials, we get $m\geq n$. Thus $n=\grk(\phi)$.\\ 

(II) Implication $i)\Longrightarrow ii)$.\\
First of all we may assume that $A$ is regular by using Lemma \ref{lemma5} and Lemma \ref{redfinite}.\\
Using Theorem \ref{structuretheorem}, we may assume that $\phi(x_i)=y_i^{p^{\a_i}}u_i$ for $1\leq i\leq \grk(\phi)$ and $\phi(x_i)=0$ if $i>\grk(\phi)$. In this case  $ii)$ is satisfied by taking $a=\max_i{p^{\a_i}}$ and $b=0$.\\

\textbf{Step 2:}
Assume that $\k=\k'$ and $\k$ is a finite field.\\
As before we may reduce to the case $A=\k[[x_1,...,\,x_n]]$ and $B=\k[[y_1,...,\,y_m]]$. Let $t$ be a variable over $\k$ and let $\K:=\k(t)$. Let $A':=\K[[x_1,...,\,x_n]]$ and $B':=\K[[y_1,...,\,y_m]]$. The homomorphism $\phi$ extends to a homomorphism $\phi_{\K} : A'\lgw B'$ in an obvious way.
According to Corollary \ref{field_tr} $\grk(\phi)=\grk(\phi_{\K})$.\\
 On the other hand let us denote by $\wdt{\phi}$ the homomorphism  $\k[[x_1,...,x_n]][\K]\lgw \k[[y_1,...,y_m]][\K]$ induced by $\phi$ (where $\k[[x_1,...,x_n]][\K]$ is the  image of $\K\otimes_{\k}\k[[x_1,...,x_n]]\lgw \K[[x_1,...,x_n]]$ defined by $\l\otimes f\lgm \l f$). Let $f\in \k[[x_1,...,x_n]][\K]$. By multiplying $f$ by an element of $\K$ we may assume that $f=\sum_{i=0}^rf_it^i$ with $f_i\in\k[[x_1,...,x_n]]$. Then $\wdt{\phi}(f)\in (y)^k$ (resp. $f\in (x)^k$) if and only if $\phi(f_i)\in(y)^k$ (resp. $f_i\in (x)^k$) for $1\leq i\leq r$  and any $k\in \N$. Thus $\phi$ satisfies $ii)$ if and only if $\wdt{\phi}$ satisfies $ii)$. It is clear that $\phi_{\K}$ satisfies $ii)$ if and only if $\wdt{\phi}$ satisfies $ii)$, $\phi_{\K}$ being the completion of $\wdt{\phi}$.\\
 Thus we use Step 1 to conclude.\\

\textbf{Step 3:}
Assume  that $\k\neq \k'$. Using Lemma \ref{lemma5} and  Cohen's Theorem (for example Corollary 31.6 of \cite{Na}), we may find an injective finite homomorphism of $\k$-algebras $\s :  A'\lgw A$ such that $\grk(\phi)=\grk(\phi\circ\s)$ and such that $A'$ is regular. By Lemma \ref{redfinite} we can replace $A$ by $A'$. So we assume that $A=\k[[x_1,...,\,x_n]]$ and $B=\k'[[y_1,...,\,y_m]]$. We denote by $A_{\k'}$ the $\k'$-algebra $A\wdh{\otimes}_{\k}\k'=\k'[[x_1,...,\,x_n]]$. We denote by $\phi_{\k'}$ the homomorphism $A_{\k'}\lgw B$ induced by $\phi$. Because $\k\lgw \k'$ is  finite, then $\grk(\phi)=\grk(\phi_{\k'})$ by Lemma \ref{lemma6}. Using Lemma 5.4 \cite{Iz3}, we see that $\phi$ satisfies $ii)$ if and only if $\phi_{\k'}$ satisfies $ii)$. Then the result follows from Step 2.$\quad\Box$\\
\\
Finally, following W. F. Osgood \cite{Os}, S. S. Abhyankar \cite{Ab} and A. Gabrielov \cite{Ga1},  we give an example  of injective homomorphisms of local rings for which the growth of the Chevalley function  is greater than any given increasing function $\a$ :
\begin{ex}\label{ex1}
Let $\a\ :\ \N\lgw \N$ be an increasing function and let $\k$ be a field. Let $(n_i)_i$ be a sequence of natural numbers such that $n_{i+1}>\a(n_i)$ for any $i$ and such that the element $\xi(Y):=\sum_{i\geq 1}Y^{n_i}$ is transcendental over $\k(Y)$ (such an element exists according to the constructive proof of Lemma 1 in \cite{ML-S}). Let us define the homomorphism $\phi : A:=\k[[x_1,\,x_2,\,x_3]]\lgw B:=\k[[y_1,\,y_2]]$ by
$$(\phi(x_1),\,\phi(x_2),\,\phi(x_3))=(y_1,\,y_1y_2,\,y_1\xi(y_2)).$$
Because $1,\,y_2,\,\xi(y_2)$ are algebraically independent over $\k$, $\phi$ is injective (cf. Part 1 of \cite{Ab}): indeed, let $f\in\Ker(\phi)$. We write $f=\sum_d f_d$, where $f_d$ is a homogeneous poynomial of degree $d$. Then $\phi(f)=\sum y_1^df_d(1,\,y_2,\,\xi(y_2))=0$. Hence, we have $f_d(1,\,y_2,\,\xi(y_2))=0$ for all $d$. This implies that $f_d=0$ for all $d$ because $1,\,y_2,\,\xi(y_2)$ are algebraically independent.\\
For any positive natural number $i$ we define:
$$f_i:=x_1^{n_i-1}x_3-\left(x_2^{n_1}x_1^{n_i-n_1}+\cdots+x_2^{n_{i-1}}x_1^{n_i-n_{i-1}}+x_2^{n_i}\right).$$
Then we get:
$$\phi(f_i)=y_1^{n_i}\xi(y_2)-y_1^{n_i}\sum_{k=1}^{i}y_2^{n_k}\in\m_B^{n_i+n_{i+1}}\subset\m_B^{\a(n_i)}$$
But $f_i\notin \m_A^{n_i+1}$ thus $\b(n_i+1)>\a(n_i)$ where $\b$ is the Chevalley function associated to $\phi$. Because $n_i\lgw+\infty$ when $i\lgw +\infty$, we get $\limsup \frac{\b(n)}{\a(n)}\geq 1$. 
 
\end{ex}

\subsection{Chevalley function and diophantine approximation}
The aim of this section is to give an interpretation in terms of diophantine approximation of the fact that the Chevalley function of a homomorphism of complete local rings is not bounded by an affine function as soon as $\phi$ is not regular.\\
Let $\phi : \k[[x_1,...,\,x_n]]\lgw \k[[y_1,...,\,y_m]]$ be a homomorphism of local $\k$-algebras, where $\k$ is an infinite field. Let us assume that $\grk(\phi)=n-1$ and that $\phi$ is injective. Using Theorem \ref{structuretheorem}, there exists a commutative diagram as follows:
$$\xymatrix{\k[[x]] \ar[r]^{\phi} \ar[d]^{\s_1}& \k[[y]] \ar[d]^{\s_2}\\
\k[[x]]\ar[r]^{\ovl{\phi}} & \k[[y]]}$$
 such that the homomorphism $\s_1$ is a composition of homomorphisms of $\k[[x]]$ defined in Definition \ref{modifications}
 and such that the homomorphism $\ovl{\phi}$ satisfies
$$\ovl{\phi}(x_i)=y_i^{p^{\a_i}}u_i \text{ for some units }u_i\in\k[[y]] \text{ and  } \a_i\in\N, \text{ for } i\leq n-1, \text{ if }\cha(\k)=p>0$$
$$\text{or }\ \ovl{\phi}(x_i)=y_i \text{ for }i\leq n-1, \text{ if } \cha(\k)=0$$ 
$$\text{ and }\ovl{\phi}(x_n)=0$$ Moreover, if $\cha(\k)=p>0$, for any $i$, $u_i=1$ whenever $\a_i=0$, and $\ini(u_i)=1$ and $u_i\notin \k[[y]]^p$ whenever $\a_i>0$.\\
The homomorphism $\s_2$ is injective and $\grk(\phi)=\grk(\s_2\circ\phi)$, thus $\grk(\s_2\circ\phi)=n-1$. From now on we will replace $\phi$ by $\s_2\circ\phi$. Hence we have the following commutative diagram:
$$\xymatrix{\k[[x]]\ar[rr]^{\phi_0:=\phi}\ar[d]^{\psi_1} &  &\k[[y]]\\
\k[[x]] \ar[rru]^{\phi_1}\ar[d]^{\psi_2} & & \\
\vdots \ar[d]^{\psi_l} &  & \\
\k[[x]]\ar[rruuu]_{\phi_l:=\ovl{\phi}} &  &}$$
where $\phi_l:=\ovl{\phi}$, and  $\psi_j$, for $1\leq j\leq l$, is one of the  homomorphisms used in the modifications  (1), (3) and (4) of Definition \ref{modifications} (resp. called homomorphisms of type (1), (3) and (4)).\\
We can remark that if $\phi_{j+1}$ is not injective and $\psi_{j+1}$ is a homomorphism of type (1) or (3), then $\phi_j$ is neither injective:\\ 
It is trivial for homomorphisms of type (1). If $\psi_{j+1}$ is a homomorphism of type (3), let $f\in\Ker(\phi_{j+1})$ and let us write $d=p^re$ with $e\wedge p=1$. Then let us define $g:=\prod_{\e\in\mathbb{U}_e}\left(f(\e x_1,\,x_2,...,\,x_d)\right)^{p^r}$ where $\mathbb{U}_e$ is the multiplicative group of the $e$-roots of unity in a finite extension of $\k$. Then $g\in\Im(\psi_{j+1})$. Let $g'\in\k[[x]]$ such that $\psi_{j+1}(g')=g$. Then $\phi_{j}(g')=\phi_{j+1}(g)=0$. Thus $\phi_j$ is not injective.\\
Nevertheless, if $\psi_{j+1}$ is a homomorphism of type (4), $\phi_j$ may be injective while $\phi_{j+1}$ is  not injective. Let us assume that $\phi_j$, for $1\leq j <k$, is injective and $\phi_{k}$ is not injective. In particular $\psi_{k}$ is a homomorphism of type (4). Because $\phi_{k}$ is not injective, we have $\dim(\k[[x]]/\Ker(\phi_{k}))=n-1=\grk(\phi_{k})$ thus there exist $a$, $b$ such that $\phi_{k}^{-1}((y)^{an+b})\subset \Ker(\phi_{k})+(x)^n$ for any $n\in\N$ according to Theorem \ref{maintheorem}.\\
Let us remark that, for any $j$, there exist $a_j\geq 1$ and $b_j\geq 0$ such that $\psi_j^{-1}((x)^{a_jn+b_j})\subset (x)^n$ for any $n\in\N$. Then there exist $a'\geq 1$ and $b'\geq 0$ such that $(\psi_{k-1}\circ\cdots\circ\psi_1)^{-1}((x)^{a'n+b'})\subset (x)^n$ for any $n\in\N$. If $\b_0$ denotes the Chevalley function of $\phi$ and $\b_{k-1}$ denotes the Chevalley function of $\phi_{k-1}$, then $\b_0(n)\leq \b_{k-1}(a'n+b)$ for any $n\in\N$. Thus $\b_0$ is not bounded by an affine function because $\b_{k-1}$ is not bounded by an affine function. We will investigate the reason why $\b_{k-1}$ is not bounded by an affine function.\\
We have to consider the following situation (here $\phi$ represents $\phi_{k-1}$, $\wdt{\phi}$ represents $\phi_k$ and $q$ represents $\psi_k$): we have the following commutative diagram
$$\xymatrix{\k[[x]] \ar[r]^{\phi} \ar[d]^{q}& \k[[y]]\\
\k[[x]]\ar[ru]^{\wdt{\phi}} & }$$
where $\phi$ is injective and $r_1(\phi)=n-1$ ; $q$ is the homomorphism  defined by $q(x_i)=x_i$ for $i\neq 1$ and $q(x_1)=x_1x_2$. Moreover $\wdt{\phi}$ is not injective: $\grk(\wdt{\phi})=n-1=\dim(\k[[x]]/\ker(\wdt{\phi}))$. From Theorem \ref{maintheorem}, there exist $a\geq 1$, $b\geq 0$ such that $\wdt{\phi}^{-1}((y)^{an+b})\subset \Ker(\wdt{\phi})+(x)^n$ for any $n\in\N$.
Let $\wdt{z}\in\k[[x]]$ be a generator of $\Ker(\wdt{\phi})$. 
 Let us denote, for any $g\in\k[[x]]$,
$$\nu_{\wdt{z}}(g):=\max\{k\in\N\ /\ g\in (\wdt{z})+(x)^k\}$$
with the assumption $\nu_{\wdt{z}}(g)=+\infty$ if $g\in(\wdt{z})$.
In particular $\nu_{\wdt{z}}(g)\leq \ord(\wdt{\phi}(g))\leq a\nu_{\wdt{z}}(g)+b$ for any $g\in \k[[x]]$.
Then $\b$ is the Chevalley function of $\phi$ means exactly  the following:
$$\begin{array}{cc} & \forall f\in\k[[x]]\  \ord(\phi(f))\leq \b(\ord(f)) \text{ and}\\
 &\forall n\in\N\ \  \exists f_n\in\k[[x]]\ /\   \ord(f_n)=n,\ \ord(\phi(f_n))=\b(\ord(f_n)).\end{array}$$
 This is equivalent to the fact  that there is a function $\g\ :\ \N\lgw \N$ such that
\begin{equation}\label{nu}\begin{array}{c}\forall  f\in\k[[x]]\ \ \nu_{\wdt{z}}(q(f))\leq \g(\ord(f)),\\
\forall n\in\N\ \  \exists f_n\in\k[[x]]\ /\   \ord(f_n)=n,\ \nu_{\wdt{z}}(f_n)=\g(\ord(f_n))\\
\text{and } \forall n\in\ N,\ \g(n)\leq \b(n)\leq a\g(n)+b.\end{array}\end{equation}\\\
Let us consider the following three rings along with the canonical injections
$$\xymatrix{A:=\k[[x_1,...,\,x_n]]\ar[r]^{i_1}  & B:=\frac{\k[[x_1,....,\,x_n]][t]}{(x_1-tx_2)}\ar[r]^{i_2} & C:= \frac{\k[[x_1,....,\,x_n,\,t]]}{(x_1-tx_2)}.}$$
The homomorphism $\t\ :\ C\lgw \k[[x_1,...,\,x_n]]$ defined by $\t(x_1)=x_1x_2$, $\t(x_i)=x_i$ for $i>1$ and $\t(t)=x_1$ is an isomorphism and $\t\circ i_2\circ i_1=q$. We will often omit the notations $i_1$ and $i_2$ in rest of the paper.\\
Let us remark the following fact:
$$\textit{The element }\t^{-1}(\wdt{z})\in C \textit{ is not algebraic over }A.$$
Indeed, if $\wdt{z}$ was  algebraic over $A$, then we would have a relation $a_0+a_1\wdt{z}+\cdots+a_d\wdt{z}^d=0$, such that $a_i\in A$ for $0\leq i\leq d$ and $a_d\neq 0$. Because $C$ is an integral domain, we amy assume that $a_0\neq 0$ by assuming that $d$ is minimal. Thus we get $\phi(a_0)=0$, because $\wdt{\phi}(\wdt{z})=0$, thus $\phi$ would not be injective which would contradicting the hypothesis. Hence $\wdt{z}$ is not algebraic over $A$.\\
\\
Let us denote by $\nu_1$ the valuation on $A$ defined by its maximal ideal  and let us denote by $\nu_2$ the valuation on $C$ defined by its maximal ideal. We still denote by $\nu_2$ its restriction on $A$ or $B$. Let us remark that $\nu_1(f)\leq \nu_2(f)\leq 2\nu_1(f)$ for any $f\in A$ and that $\nu_2(f)=\nu_1(q(f))$ for any $f\in A$.\\
Let us denote by $\K_A$ (resp. $\K_C$) the field of fractions of $A$ (resp.  $C$). Let us remark that $\K_A$ is also the field of fractions of $B$.
Let us denote, for any $f\in \K_A$, $|f|_1:=e^{-\nu_1(f)}$, and for any $g\in \K_C$ let us denote $|g|_2:=e^{-\nu_2(g)}$. Then $|.|_1$  and $|.|_2$ are non-archimedian norms on $\K_A$ and $\K_C$ respectively. Let us denote by $\wdh{\K}$ the completion of $\K_A$ with respect to $|.|_2$. We can remark that there is a natural injection $\K_C\hookrightarrow \wdh{\K}$.\\
Let us come back to $\wdt{z}$, the generator of $\Ker(\wdt{\phi})$. We have the following lemma:
\begin{lemma}
The element $\t^{-1}(\wdt{z})$ satisfies the following property: there exists a decreasing function $\a\ :\ \R^+\lgw \R^+$ such that
$$\left|\frac{f}{g}-\t^{-1}(\wdt{z})\right|_2\geq \a(|g|_2) \ \ \forall f\in A,\, g\in B.$$
Moreover, the Chevalley function of $\phi$ is not bounded by an affine function because of the following fact: If $\a$ is the greatest function satisfying the above inequality, then $\frac{\ln(\a(u))}{\ln(u)}\lgw 0$ as $u$ goes to 0.
\end{lemma}
\begin{proof}The fact (\ref{nu}) means that for any $f\in A$ and for any  $g\in C$, we have $\nu_2(f-g\t^{-1}(\wdt{z}))\leq \g(\nu_1(f))$ and this inequality is the best possible. This is equivalent to $\nu_2(f-g\t^{-1}(\wdt{z}))\leq \g(\nu_1(f))$ for any $f\in A$ and any $g\in B$, where the inequality is the best possible, because $C$ is the completion of $B$ for $\nu_2$. Thus for any $f\in A$ and for any  $g\in B$, we have \begin{equation}\label{ine}\nu_2(f-g\t^{-1}(\wdt{z}))\leq \g'(\nu_2(f)),\end{equation} with $\g(\frac{n}{2})\leq \g'(n)\leq \g(n)$ for any $n\in\N$, and this inequality is the best possible. We do not make any restriction if we assume that $\nu_2(f)=\nu_2(g\t^{-1}(\wdt{z}))$: if it is not the case we have $\nu_2(f-g\t^{-1}(\wdt{z}))\leq \nu_2(f)$, but clearly the least function $\g'$ satisfying the inequality (\ref{ine}) for any $f$ and $g$ satisfies $\g'(n)\geq n$ for any $n\in\N$. Thus we get
$\left|\frac{f}{g}-\t^{-1}(\wdt{z})\right|_2\geq \a(|g|_2)$ for any $f\in A,\,g\in B$ with $\a(u):=e^{-\g'(\ln(u))}$ for any $u>0$. We get $\frac{\ln(\a(u))}{\ln(u)}\lgw 0$ as $u$ goes to 0, because $\g'$ is not bounded by an affine function, this following from the fact that $\b$, thus $\g$, is neither bounded by an affine function.
\end{proof}

\begin{rmk}
Let us remark the following fact: if $z\in C$ is algebraic over $A$, then there does not exist any function $\a\ :\ \R^+\lgw \R^+$ such that $\left|\frac{f}{g}-z\right|_2\geq \a(|g|_2) \ \ \forall f\in A,\, g\in B.$ Indeed since $z$ is algebraic over $A$ there exists a relation 
$a_dz^d+\cdots a_1z+a_0$ such that $a_i\in A$ for $0\leq i\leq d$ and $a_d\neq 0$. Because $C$ is an integral domain, we may assume that $a_0\neq 0$. Thus $z.w=a_0$ with $w:=-(a_dz^{d-1}+\cdots+a_1)\in C$. For any $n\in \N$, let us denote by $w_n$ an element of $B$ such that $w_n-w\in \m_C^n$ and $\nu_2(w_n)=\nu_2(w)$. Such a $w_n$ exists because $C$ is the completion of $B$. Thus we have $\nu_2\left(z-\frac{a_0}{w_n}\right)=n-\nu_2(w)$ for any $n\in N$. Thus $\left|z-\frac{a_0}{w_n}\right|_2\lgw 0$ as $n\lgw \infty$, but $|w_n|_2=|w|_2\neq 0$ for any $n\in N$.
\end{rmk}

\section{Homomorphisms of Henselian  $\k$-algebras}\label{section-Henselian}
In this section and the next one, we  study a particular example of homomorphisms of local $\k$-algebras: namely the homomorphisms of W-system. Such homomorphisms generalize homomorphisms of analytic local rings in the sense that the local rings that we consider satisfy the Weierstrass Division Theorem. In particular we have been inspired by the work of S. S. Abhyankar and M. van der Put \cite{A-vdP} on analytic $\k$-algebras.\\
\subsection{Terminology}
From now on we assume that $\k$ is a field  of any characteristic. 
\begin{definition}\label{W-sys}\cite{D-L} By a \emph{Weierstrass System of local $\k$-algebras}, or a \emph{W-system} over $\k$, we mean a family of  $\k$-algebras $\k\llceil x_1,...,\,x_n\rrceil$, $n\in\N$ such that:
\begin{enumerate}
\item[i)] For $n=0$, the $\k$-algebra is $\k$,\\
For $n\geq 1$, $\k[x_1,...,\,x_n]_{(x_1,...,\,x_n)}\subset \k\llceil x_1,...,\,x_n\rrceil \subset \k[[x_1,...,\,x_n]]$\\
and $\k\llceil x_1,...,\,x_{n+m}\rrceil \cap \k[[x_1,...,\,x_n]]=\k\llceil x_1,...,\,x_n\rrceil$ for $m\in\N$.
For any permutation of $\{1,...,\,n\}$, denoted by $\s$, $\k\llceil x_{\s(1)},...,\,x_{\s(n)}\rrceil=\k\llceil x_1,...,\,x_n\rrceil$.
\item[ii)]  Any element of $\k\llceil x\rrceil$, $x=(x_1,...,\,x_n)$, which is a unit in $\k[[x]]$, is a unit in $\k\llceil x\rrceil$.
\item[iii)] Let $f\in(x)\k\llceil x\rrceil$ such that $f(0,...,\,0,\,x_n)\neq 0$. We denote $d:=\ord_{x_n}f(0,...,\,0,\,x_n)$. Then for any $g\in\k\llceil x\rrceil$ there exist a unique $q\in\k\llceil x\rrceil$ and a unique $r\in\k\llceil x_1,...,\,x_{n-1}\rrceil[x_n]$ with $\deg_{x_n}r<d$ such that $g=qf+r$.
\item[iv)] (if $char(\k)>0$) If $h\in(y_1,...,\,y_m)\k[[y_1,...,\,y_m]]$ and $f\in\k\llceil x_1,...,\,x_n\rrceil$ such that $f\neq 0$ and $f(h)=0$, then there exists $g\in\k\llceil x\rrceil$ irreducible in $\k\llceil x\rrceil$ such that $g(h)=0$ and such that there does not exist any unit $u(x)\in\k\llceil x\rrceil$ with $u(x)g(x)=\sum_{\a\in\N^n}a_{\a}x^{p\a}$ ($a_{\a}\in\k$).
\end{enumerate}
\end{definition}
\begin{rmk}\label{rmk-W} Let $\k\llceil x\rrceil$ be a W-system over $\k$.
\begin{enumerate}
\item[i)] From \cite{D-L} the ring $\k\llceil x_1,...,,x_n \rrceil$ ($n\in\N$) is a Noetherian regular local ring with maximal ideal $(x_1,...,,x_n)$ and its completion at its maximal ideal  is $\k[[x]]$.
\item[ii)] For any $f\in\k\llceil x_1,...,x_{n+m}\rrceil$ and any $g_1,...,\,g_m\in(x)\k\llceil x_1,...,\,x_n\rrceil$, $$f(x_1,...,\,x_n,\,g_1(x),...,,g_m(x))\in\k\llceil x_1,...,,x_n\rrceil.$$ (\cite{D-L})
\item[iii)] For any $f\in\k\llceil x\rrceil$, if there is $g\in \k[[x]]$ such that $f=x_1g$, then $g\in\k\llceil x\rrceil$. (\cite{D-L})
\item[iv)] From Theorem 44.4 \cite{Na}, iii) implies that $\k\llceil x\rrceil$ is a Henselian local ring. In fact it is proven in \cite{D-L} that $k\llceil x\rrceil$ has the Artin Approximation Property, and by \cite{Po} and \cite{Ro} (where it is proven that a local Noetherian ring has the Artin Approximation Property if and only if it is Henselian and excellent), we see that $\k\llceil x\rrceil$ is excellent. In \cite{D-L}  (Remark 10) it is said that if a family of excellent rings satisfies i), ii) and iii), then it satisfies iv).
\item[v)] Let $d>1$, $d\wedge \cha(\k)=1$,  let $a\in\k^*$ be a $d$-th power in $\k$ and let $f(x)\in(x)\k\llceil x\rrceil$. It means that $P(T)=T^d-(a+f(x))\in\k\llceil x\rrceil[T]$ has a non-zero solution modulo $(x)$. Thus, because $\k\llceil x\rrceil$ is Henselian, $P(T)$ has a solution in $\k\llceil x\rrceil$. Hence $a+f(x)$ has a $d$-th root in $\k\llceil x\rrceil$.\\
If $d=\cha(\k)>0$, $a\in\k^*$ is a $d$-th power in $\k$ and  $f(x)\in(x)\k\llceil x\rrceil$ is a $d-$power in $\k[[x]]$, then $P(T)=T^d-(a+f(x))$ has $d$-root in $\k[[x]]$, thus it has a $d$-root in $\k\llceil x\rrceil$ by the Artin Approximation Theorem \cite{D-L}.\end{enumerate}
\end{rmk}

In fact we can give a quick proof of the fact that W-systems satisfy the Artin Approximation Property if we assume that the rings of the family are excellent, using the Popescu's smoothing Theorem (cf. \cite{Po} or \cite{Sp2}):
\begin{theorem}\cite{D-L}
Let $\k\llceil x\rrceil$ be a W-system over $\k$ and let us assume that $\k\llceil x_1,...,\,x_n\rrceil$ is excellent for any $n\in\N$. Then for any $f=(f_1,...,\,f_p)\in\k\llceil x,\,y\rrceil$ with $x=(x_1,...,\,x_n)$ and $y=(y_1,....,\,y_m)$, for any $c\in\N$ and for any $\ovl{y}\in(x)\k[[x]]^m$ such that $f(\ovl{y})=0$ there exist $y_c\in(x)\k\llceil x\rrceil^m$ such that $f(y_c)=0$ and $\ovl{y}_i-y_{c,i}\in (x)^c$.
\end{theorem}
\begin{proof}
 We may assume that $p=1$ by replacing $(f_1,...,\,f_p)$ by 
 $$f:=f_1^2+x_1(f_2^2+x_1(f_3^2+x_1(\cdots +x_1f_p^2)^2)^2\cdots)^2.$$ 
 By assumption there exist $\ovl{h}_{i}(x,\,y)\in\k[[x,\,y]]$, $1\leq i\leq n$, such that
$$f(y)+\sum_{i=1}^n(y_i-\ovl{y}_i(x))\ovl{h}_i(x,\,y)=0.$$
Because the ring $\k\llceil y\rrceil\langle x\rangle$ is Henselian and excellent, it satisfies the Artin Approximation Property for algebraic equations \cite{Po} ($\k\llceil y\rrceil\langle x\rangle$ is the Henselization  of $\k\llceil y \rrceil[x]_{(y, x)}$). Thus there exist $f_i(x,\,y),\,h_i(x,\,y)\in \k\llceil y\rrceil\langle x\rangle^n$, $1\leq i\leq n$, such that
$$\ovl{h}_i(x,\,y)-h_i(x,\,y),\,f_i(x,\,y)-\ovl{y}_i(x)\in (x,\,y)^c, 1\leq i\leq n$$
$$\text{ and }f(y)+\sum_{i=1}^n(y_i-f_i(x,\,y))h_i(x,\,y)=0.$$
We may assume that $c\geq 2$. In this case the Jacobian matrix of $(y_i-f_i(x,\,y),\ 1\leq i\leq n)$ with respect to $y_1,...,\,y_n$ has determinant equal to 1 modulo $(x,\,y)$. The Henselian property asserts that there exist $y_{i,c}(x)\in\k\llceil x \rrceil$ such that
$$y_{i,c}(x)-f_i(x,\,y_{1,c}(x),...,\,y_{n,c}(x))=0\text{ for } 1\leq i\leq n.$$
Then
$$f(y_{1,c}(x),...,\,y_{n,c})=0$$
$$\text{ and } y_{i,c}(x)-\ovl{y}_i\in (x)^c,\ 1\leq i\leq n.$$
\end{proof}

\begin{rmk}
In the same way we may prove that the rings $\k\llceil x\rrceil$ satisfy the Strong Artin Approximation Property (cf. Theorem 7.1 \cite{D-L}) using the fact that a ring that satisfies the Artin Approximation Property satisfies also the Strong Artin Approximation Property  \cite{P-P}.
\end{rmk}

\begin{ex}\label{ex_W}
\begin{enumerate}
\item[i)] The family $\k[[x_1,...,\,x_n]]$ is a W-system over $\k$.
\item[ii)] Let $\k\langle x_1,...,\,x_n\rangle$ be the Henselization of the localization of $\k[x_1,...,\,x_n]$ at the maximal ideal $(x_1,...,\,x_n)$. Then, for $n\geq 0$, the family $\k\langle x_1,...,\,x_n\rangle$ is a W-system over $\k$.
\item[iii)] The family $\k\{x_1,...,\,x_n\}$ (the ring of convergent power series in $n$ variables over a valued field $\k$) is a W-system over $\k$.
\item[iv)] The family of Gevrey power series in $n$ variables over a valued field is a W-system over $\k$  \cite{Br}.

\end{enumerate}
\end{ex}

\begin{definition}
Let $\k\llceil x\rrceil$ be a W-system over a field $\k$. A local ring $A$ is a \emph{local $\k$-algebra with respect to this W-system}  if $A$ is isomorphic to $\k\llceil x_1,...,\,x_n\rrceil[\k']/I$ for some $n>0$, where $\k'$ is a finite field extension of $\k$  and $I$ is an ideal of $\k\llceil x_1,...,\,x_n\rrceil[\k']$ ($\k\llceil x\rrceil[\k']$ is the image of the $\k$-homomorphism $\k\llceil x \rrceil[t_1,...,,t_s]\lgw \k'[[x]]$ where $x_i$ is sent on $x_i$ and $t_j$ is sent on $\e_j$, where $\e_1$,..., $\e_s$ is a $\k$-basis of $\k'$).\\
A homomorphism of local $\k$-algebras $A\lgw B$ is called a \emph{homomorphism of Henselian $\k$-algebras} if $A$ and $B$ are local $\k$-algebras with respect to the same W-system over $\k$ and the homomorphism is a homomorphism of local $\k$-algebras.
\end{definition}
\begin{rmk}
If $A$ is a local $\k$-algebra with respect to a W-system, then its residue field  is a finite extension of $\k$. If $A\lgw B$ is a homomorphism of Henselian $\k$-algebras, then the residue field of $B$ is a finite extension of the residue field of $A$.
\end{rmk}
 \begin{rmk} Since an integral extension of a local Henselian ring is an Henselian ring (\cite{Na} 43.16), any  local $\k$-algebra with respect to a W-system is a  Henselian ring.\\
 Let $\k\llceil x\rrceil$ be a W-system.  Thus, from \cite{La}, $\k\llceil x \rrceil[\k'] $ satisfies Property iii) of the definition of a W-system if $\k\lgw \k'$ is a finite field  extension. Moreover it is straightforward to show that $\k\llceil x\rrceil[\k']$ satisfies i) and ii) in the definition of a W-system. Finally, from Remark \ref{rmk-W} iv), and since any finite extension of an excellent ring is an excellent ring, we see that $\k\llceil x\rrceil[\k']$ satisfies iv) of the definition of a W-system. Hence $\k\llceil x\rrceil[\k']$ is a W-system with respect to $\k'$ if $\k\llceil x \rrceil$ is a W-system over $\k$ and $\k\lgw \k'$ is a finite field extension. 
\end{rmk}
\begin{definition}A  homomorphism $\phi : A\lgw B$ of Henselian  $\k$-algebras is $\emph{strongly injective}$ if the map $\wdh{A}/A\lgw \wdh{B}/B$ induced by $\phi$ is injective (or equivalently if $\wdh{\phi}^{-1}(B)=A$).
\end{definition}
Finally we give the following version of the Weierstrass Preparation Theorem:
\begin{prop}\label{Weierstrass}\textit{Weierstrass Preparation Theorem}\\
Let $A$ and $B$ be local $\k$-algebras with respect to a W-system  and $\phi : A\lgw B$ be a homomorphism of Henselian  $\k$-algebras. Then $\phi$ is finite if and only if $\phi$ is quasi-finite (i.e. $B/\m_A B$ is finite over $A/\m_A$). 
\end{prop}

\begin{proof}
It is well known that ii) of Definition \ref{W-sys} is equivalent to the proposition when there exist surjective homomorphisms $\k\llceil x_1,...,\,x_n\rrceil\lgw A$ and $\k\llceil y_1,...,\,y_m\rrceil\lgw B$ for some W-system $\k\llceil x\rrceil$ (\cite{To} or \cite{Ab2} for example). Because $\k\llceil x\rrceil[\k']$ is a W-system over $\k'$ as soon as $\k\llceil x\rrceil$ is a W-system over $\k$ and $\k\lgw \k'$ is a finite extension of fields, the proposition is proven.
\end{proof}

\begin{corollary}\label{regular}
Let $A$ be a regular  local $\k$-algebra with respect to a W-system $\k\llceil x\rrceil$ and let $(a_1,...,a_n)$ be a regular system of parameters of $A$. Let $\k'$ be the coefficient field of $A$.  Let $\phi : \k\llceil x_1,...,\,x_n\rrceil[\k']\lgw A$  be the unique homomorphism of local $\k'$-algebras such that $\phi(x_i)=a_i$ for $1\leq i\leq n$. Then $\phi$ is an isomorphism.
\end{corollary}

\begin{proof}
Follows from Proposition \ref{Weierstrass}.
\end{proof}

\subsection{Strongly injective homomorphisms}
We state now the following results about homomorphisms of Henselian $\k$-algebras:

\begin{lemma}\label{finite-injective}(\cite{A-vdP} Lemma 2.1.2)
Let $\phi$ be a homomorphism of Henselian  $\k$-algebras. If $\phi$ is injective and finite, then $\wdh{\phi}$ is injective and finite and $\phi$ is strongly injective. 
\end{lemma}
\begin{proof}
Local  $\k$-algebras are Zariski rings (cf. Theorem 9, chapter VIII of \cite{Z-S}). Then, using Theorem 5 and Theorem 11 Chapter VIII of \cite{Z-S}, we see that $\wdh{\phi}$ is finite and injective. Then using Theorem 15 Chapter VIII of \cite{Z-S}, we see that $\phi$ is strongly injective.
\end{proof}

Let $\k$ be a field of any characteristic and $\k\llceil x \rrceil$ be a W-system with respect to $\k$. We define the local $\k$-homomorphism $q\ :\ \k\llceil x_1,...,\,x_n\rrceil\lgw \k\llceil x_1,...,,x_n\rrceil$ by $q(x_1)=x_1x_2$  and $q(x_i)=x_i$ for $i>1$. For any $d\in\N^*$ we define the local $\k$-homomorphism $\psi_d\ :\ \k\llceil x_1,...,,x_n\rrceil \lgw \k\llceil x_1,...,\,x_n\rrceil$ by $\psi_d(x_1)=x_1^d$ and $\psi_d(x_i)=x_i$ for $i\neq 1$.

\begin{lemma}\label{st-inj_quad}
Let $A$ be a $\k$-algebra with respect to a W-system  denoted by $\k\llceil x_1,...,\,x_n\rrceil$. Any composition of $\k$-automorphisms of $A$ and of homomorphisms of the form $q$ and $\psi_d$ is injective. Moreover  $\k$-automorphisms and homomorphisms of the form $\psi_d$ are strongly injective.
\end{lemma}

\begin{proof}
It is clear that $\k$-automorphisms,  homomorphisms $\psi_d$ and $q$ are injective. Moreover it is clear that $\k$-automorphisms are strongly injective.\\
Let $f(x_1,...,\,x_n)\in\k[[x_1,...,\,x_n]]$ such that $\psi_d(f)=f(x_1^d,\,x_2,...,,x_n)\in\k\llceil x \rrceil$. We have
$$f(x_1^d,\,x_2,...,\,x_n)=(x_1^d-y)q(x,\,y)+r(x,\,y)$$
with $r(x,y)\in\k\llceil x_2,...,\,x_n,\,y\rrceil[x_1]$ and $\deg_{x_1}(r)<d$ according to iii) of Definition \ref{W-sys}. On the other hand, 
$$f(x_1,...,\,x_n)=(x_1-y)\ovl{q}(x,\,y)+\ovl{r}$$
with $\ovl{r}\in\k[[x_2,...,\,x_n,\,y]]$ according to the  formal Weierstrass division Theorem. Thus
$$f(x_1^d,\,x_2,...,\,x_n)=(x_1^d-y)\ovl{q}(x_1^d,\,x_2,...,\,x_n,\,y)+\ovl{r}$$
and because the division is unique in $\k[[x,\,y]]$, we see that $q(x,\,y)=\ovl{q}(x_1^d,\,x_2,...,\,x_n,\,y)$ and $r=\ovl{r}$ does not depend on $x_1$. Since $r(x_2,....,\,x_n,\,x_1)=f(x_1,...,,x_n)$, $f\in\k\llceil x\rrceil$.
\end{proof}

\begin{definition}\label{good_hens}
Let $\k\llceil x\rrceil$ be a W-system over $\k$. We say that it is a \emph{good W-system} if the homomorphism $q\ :\ \k\llceil x_1,...,\,x_n\rrceil\lgw \k\llceil x_1,...,,x_n\rrceil$ defined by $q(x_1)=x_1x_2$  and $q(x_i)=x_i$ for $i>1$ is strongly injective. A homomorphism of local $\k$-algebras $A\lgw B$ is called a  \emph{homomorphism of good Henselian $\k$-algebras} if $A$ and $B$ are local $\k$-algebras with respect to some good W-system and the homomorphism is local.
\end{definition}

\begin{rmk}
It would be interesting to know if any W-system is a good W-system.
\end{rmk}

\begin{lemma}
The W-systems $\k\llceil x\rrceil$ presented in Example \ref{ex_W} are good W-systems over $\k$. The same is true for the respective W-systems over $\k'$,  $\k\llceil x\rrceil[\k']$, where $\k\lgw \k'$ is a finite field extension.
\end{lemma}

\begin{proof}
It is clear for $\k[[x]]$. For the convergent power series, we have just to remark that for any $f=\sum_{\a}a_{\a}x_1^{\a_1}...x_n^{\a_n}\in\k[[x]]$ such that $q(f)$ is convergent, there exist $R_1>0$,..., $R_n>0$ such that $ \sum_{\a}|a_{\a}|R_1^{\a_1}R_2^{\a_1+\a_2}R_3^{\a_3}...R_n^{\a_n}<+\infty$. Thus $f$ is convergent because $\sum_{\a}|a_{\a}|(R_1R_2)^{\a_1}R_2^{\a_2}R_3^{\a_3}...R_n^{\a_n}<+\infty$. The proof is the same when we take a finite extension of the residue field (see \cite{A-vdP} Lemma 2.2.1). The proof is similar for Gevrey power series.\\
Let $f\in\k'[[x_1,...,\,x_n]]$. We have $f=\sum f_l\e_l$ where $f_l\in \k[[ x_1,...,\,x_n]]$ for $1\leq l\leq r$. Assume that $q(f)\in \k\langle x_1,...,\,x_n\rangle[\k']$. It is clear that $q(f_l) \in\k\langle x_1,...,\,x_n\rangle$ for $1\leq l\leq r$. It is enough to prove that if $f\in\k[[ x_1,...,\,x_n]]$ satisfies $q(f)\in \k\langle x_1,...,\,x_n\rangle$ then $f\in \k\langle x_1,...,\,x_n\rangle$. So let $f\in\k[[ x_1,...,\,x_n]]$ such that
$$g:=f(x_1x_2,\;x_2,...,\,x_n)\in\k\langle x_1,...,\,x_n\rangle.$$
There exist $s\in\N$ and $a_i\in\k[x_1,...,\,x_n]_{(x)}$ for $0\leq i\leq s$ such that
\begin{equation}\label{int2}a_sg^s+\cdots+a_1g+a_0=0.\end{equation}
We write $a_i=\sum_{\a} a_{i,\,\a}x_1^{\a_1}...x_n^{\a_n}$ for any $i$ with $a_{i,\,\a}\in \k$.  Multiplying Relation (\ref{int2}) by some power of $x_2$, we may assume that any $\a\in\N^n$ such that $a_{i,\,\a}\neq 0$ satisfies $\a_2\geq \a_1$. Then there exist  $b_i\in\k[x_1,...,\,x_n]_{(x)}$ such that $q(b_i)=a_i$ for $0\leq i\leq s$. We have $b_sf^s+\cdots+b_1f+b_0=0$, hence $f\in \k\langle x_1,...,\,x_n\rangle$ and $q$ is strongly  injective.
\end{proof}

\begin{lemma}\label{comp-strongly}(\cite{A-vdP} Lemma 2.1.3)
Let $\phi :  A\lgw B$ and $\phi' : B\lgw C$ be homomorphisms of  Henselian $\k$-algebras. If $\phi'\circ\phi$ is strongly injective then $\phi$ is strongly injective. 
\end{lemma}
\begin{proof}
Follows from the definitions.
\end{proof}

\begin{lemma}\label{analytic-grk}
Let $\k\llceil x\rrceil$ be a W-system over $\k$.
Let $\phi : \k\llceil x_1,...,\,x_n\rrceil\lgw \k\llceil y_1,...,\,y_m\rrceil[\k']$ where $\k\lgw\k'$ is a finite field extension. Let $\phi_{\k'}$ denote the induced homomorphism of Henselian $\k'$-algebras : $\k\llceil x_1,...,\,x_n\rrceil[\k']\lgw \k\llceil y_1,...,\,y_m\rrceil[\k']$. Then $\grk(\phi)=\grk(\phi_{\k'})$.
\end{lemma}

\begin{proof}
By Lemma \ref{lemma3}, we amy replace the Henselian  algebras by their completions. Then the result comes from the fact that $\k[[x_1,...,\,x_n]]\lgw \k'[[x_1,...,\,x_n]]$ is finite and Lemma \ref{lemma5}.
\end{proof}

Theorem \ref{structuretheorem} is still valid for homomorphisms of Henselian  $\k$-algebras:

\begin{theorem}\label{structuretheorem'}
Let $\k$ be an infinite field of any characteristic  and let $\k\llceil x\rrceil$ be a W-system over $\k$. Let us consider an homomorphism $\phi : A\lgw B$, where $A=\k\llceil x_1,...,\,x_n\rrceil$ and $B=\k\llceil y_1,...,\,y_m\rrceil$.  Then there exists an admissible finite sequence of homomorphisms  $ (\phi_i)_{i=0}^k : \k\llceil x_1,...,\,x_n\rrceil\lgw\k\llceil y_1,...,\,y_m\rrceil$  such that $\phi_0=\phi$. The last homomorphism $\phi_k$ satisfies
$$\phi_k(x_i)=y_i^{p^{\a_i}}u_i\ \text{ for some units }u_i\text{ if }\cha(\k)=p>0$$
$$\text{or }\ \phi_k(x_i)=y_i\ \text{ if } \cha(\k)=0,\text{ for }i\leq \grk(\phi)$$ 
$$\text{ and }\phi_k(x_i)=0\ \text{ for }i>\grk(\phi).$$ Moreover, if $\cha(\k)=p>0$, for any $i$, $u_i=1$ whenever $\a_i=0$, and $\ini(u_i)=1$ and $u_i\notin B^p$ whenever $\a_i>0$.
\end{theorem}

\begin{proof}
Modifications of types (2) and (4) are allowed according to i) of the definition of W-systems and Remark \ref{rmk-W} iii).
Steps 0, 1, 4 and 5 involve only $\k$-automorphisms of $\k[[x_1,...,\,x_n]]$ and of $\k[[y_1,...,\,y_n]]$ that are defined by polynomials. For Steps 2, 4 and 5, using modifications of type $(3)$, we  take $d$-roots of elements of $\k\llceil x\rrceil$ in $\k[[x]]$ and they are in $\k\llceil x\rrceil$ from Remark \ref{rmk-W} v). The only problem may occur at Step 3, where we replace $x_j$ by an element of the form $x'_j:=x_j-\sum_{\underline{k}}c_{\underline{k}}x_1^{k_1}...x_{j-1}^{k_{j-1}}$ such that $\phi(x'_j)=0$ (and the sum is not finite) because we do not know if $x'_j\in\k\llceil x\rrceil$. When $\cha(\k)=0$ this is obvious because $\phi(x_i)=y_i$ for $1\leq i\leq j-1$ by assumption (see Remark \ref{char0}).\\
From now on we assume that $\cha(\k)=p>0$. We assume that $A=\k\llceil x_1,...,\,x_n\rrceil$ and $B=\k\llceil y_1,...,\,y_m\rrceil$ and we will prove that $x'_j\in A$. We will use the following lemma:
\begin{lemma}\label{str-inj-struct}
Assume that $\cha(\k)=p>0$. Let us consider $\phi : \k\llceil x_1,...,\,x_{j-1}\rrceil\lgw \k\llceil y_1,...,\,y_m\rrceil$  such that   $\phi(x_i)=y_i^{p^{\a_i}}u_i$ for some units $u_i$, for $1\leq i\leq j-1$. Then $\phi$ is strongly injective.
\end{lemma}
In particular, because $\phi(x_j)=\phi\left(\sum_{\underline{k}}c_{\underline{k}}x_1^{k_1}...x_{j-1}^{k_{j-1}}\right)\in \k\llceil y_1,...,\,y_m\rrceil$,  we see that 
$$\sum_{\underline{k}}c_{\underline{k}}x_1^{k_1}...x_{j-1}^{k_{j-1}}\in\k\llceil y_1,...,\,y_m\rrceil,$$
and so $x'_j\in\k\llceil y_1,...,\,y_m\rrceil$.\end{proof}
Now we give the proof of  Lemma \ref{str-inj-struct}:

\begin{proof}[Proof of Lemma \ref{str-inj-struct}]
Let us denote by $\pi$ the quotient homomorphism $\k\llceil y_1,...,\,y_m\rrceil\lgw \k\llceil y_1,...,\,y_{j-1}\rrceil$. Then the homomorphism induced by $\pi\circ\phi$: 
$$\k\llceil x_1,...,\,x_{j-1}\rrceil/(x_1,...,\,x_{j-1})\lgw\k\llceil y_1,...,\,y_{j-1}\rrceil/(\pi\circ\phi((x_1,...,\,x_{j-1})))$$
 is finite. 
 Using Proposition \ref{Weierstrass}, we see that $\pi\circ\phi$ is finite.  Moreover, $\pi\circ\phi$ is injective because $\grk(\pi\circ\phi)=j-1$. Using Lemma \ref{finite-injective}, we see that $\pi\circ\phi$ is strongly injective, and from Lemma \ref{comp-strongly} we see that $\phi$ is strongly injective.
\end{proof}

In particular we get the following result, which is a weak version of a theorem of A. M. Gabrielov \cite{Ga2}:

\begin{theorem}\label{grk-strongly}
Let $\k$ be a field of any characteristic. Let $\phi\ :\ A\lgw B$ be a homomorphism of good Henselian $\k$-algebras, where $A$  and $B$ are regular. If $\grk(\phi)=\dim (A$), then $\phi$ is strongly injective.
\end{theorem}

\begin{proof}
Using Corollary \ref{regular}, we may assume that $A=\k\llceil x_1,...,\,x_n\rrceil$ and $B=\k\llceil y_1,...,\,y_m\rrceil[\k']$ where $\k\lgw \k'$ is finite and $\k\llceil x\rrceil$ is a W-system over $\k$. If we replace $A$ by $\k\llceil x\rrceil[\k']$ then the geometric rank will not change by Lemma \ref{analytic-grk}. Moreover if the induced homomorphism $\phi_{\k'}\ :\ \k\llceil x_1,...,\,x_n\rrceil [\k']\lgw \k\llceil y_1,...,\,y_m\rrceil[\k']$ is strongly injective, then $\phi$ is strongly injective. So from now on we assume that $\k=\k'$.\\
From Corollary \ref{field_tr} we may assume that $\k$ is an infinite field.
Using Theorem \ref{structuretheorem'}, we see that $\s_1\circ\phi=\ovl{\phi}\circ\s_2$ where the homomorphisms $\s_1$ and $\s_2$ are compositions of $\k$-automorphisms of $A$  and $B$, of homomorphisms of the form $q$ and $\psi_d$, and $\ovl{\phi}$ is defined by $\ovl{\phi}(x_i)=y_i^{p^{\a_i}}u_i$, for some units $u_i$ and some $\a_i\in\N$, for all $i$. Then using Lemmas  \ref{st-inj_quad}, \ref{comp-strongly} and \ref{str-inj-struct}, we see that $\phi$ is strongly injective.
\end{proof}

\section{Two particular cases}\label{section-ex}
\subsection{The two-dimensional case}
Example \ref{ex1} shows  that we can construct injective homomorphisms $\phi : A\lgw B$ with $\grk\phi<\dim (A)$ as soon as $\dim (A)\geq 3$. We prove here that it is not possible to find such examples when $A$ is a Henselian $\k$-algebra and $\dim (A)\leq 2$.\\
In fact, it is obvious that if  $\dim (A)=1$ and $\phi$ is injective then $\grk(\phi)=1$. Indeed,  using  Lemma \ref{lemma5}, we can replace $A$ by $\k\llceil x\rrceil$, where $x$ is a single variable and $B$ by $\k\llceil y_1,...,\,y_m\rrceil[\k']$. Then the result is immediate.\\
When $\dim (A)=2$ we have the following result that shows us that $\dim (A)=2$ is a nice case as remarked by S. S. Abhyankar and M. van der Put in \cite{Ab} and \cite{A-vdP}:

\begin{theorem}\label{dim2}
Let $\phi\ : A\lgw B$ be a homomorphism of Henselian $\k$-algebras where $\wdh{A}$ is an integral domain of dimension 2 and $B$ is regular. 
Then $\phi$ is injective if and only if $\grk(\phi)=2$.
\end{theorem}

\begin{proof}
From Lemma \ref{lemma4}, we see that $\grk(\phi)=2$ implies that $\phi$ is injective. So from now on we assume that $\phi$ is injective.\\
By Theorem 2.1 \cite{D-L} there exists an injective and finite homomorphism of Henselian $\k$-algebras
$\pi\ :\ \k\llceil x_1,\,x_2\rrceil \lgw A$ (where $\k\llceil x\rrceil$ is a W-system over $\k$), so using Lemma \ref{lemma5}, we can replace $A$ by $\k\llceil x_1,\,x_2\rrceil$. Because $B$ is regular we assume that $B=\k\llceil y_1,...,\,y_m\rrceil[\k']$ (Corollary \ref{regular}) where $\k'$ is a finite field extension of $\k$. Then, we can replace $\k\llceil x_1,\,x_2\rrceil$ by $\k\llceil x_1,\,x_2\rrceil[\k']$  using Lemma \ref{analytic-grk}.\\
Let $t$ be a variable over $\k$ and let $\K:=\k(t)$. Let $\phi_{\K}: \K[[x_1,x_2]]\lgw \K[[y_1,...,y_m]]$ be the homomorphism induced by $\phi$. If $\phi$ is injective then $\phi_{\K}$ is also injective: otherwise there would exist a sequence $(f_n)_n\in\k[[x_1,x_2]][\K]^{\N}$ such that $\phi_{\K}(f_n)\in ( y)^n$ and $f_n-f_{n+1}\in(x)^n$ for any $n\in\N$. Let $d:=\ord(f_n)$ for $n$ large enough and let us denote by $\b$ the Chevalley function of $\phi$. Let $N\in\N$ such that $\ord(\phi_{\K}(f_N))>\b(d)$. We may assume that $g:=f_N\in\k[[x_1,x_2]][\k(t)]$ by multiplying it by an element of $\K$. We write $g=\sum_{j=1}^rg_{j}t^j$ with $g_{j}\in\k[[x_1,x_2]]$ for  $0\leq j\leq r$.  Then $\phi_{\K}(g)=\sum_{j=1}^r\phi(g_{j})t^j\in (y)^{\b(d)+1}$ by assumption thus $\ord(\phi(g_j))\geq \b(d)+1$ for $0\leq j\leq r$, hence $\ord(g_j)\geq d+1$ by definition of $\b$. This contradicts $\ord(g)=d$. Hence $\phi_{\K}$ is injective and we may assume that $\k$ is infinite from Corollary \ref{field_tr}.\\
To compute $\grk(\phi)$ we use the algorithmic proof of Theorem \ref{structuretheorem}.\\
We first give the proof when $\cha(\k)=p>0$. Using Step 1, we may assume that $\phi(x_1)=y_1^{d}u$ for some unit $u$. Then we define $$\Delta:=\{\a\in \N\ /\ \exists z\in \k\llceil x_1,\,x_2\rrceil[\k']\text{ with }\ord(z)=1, \text{ and  } \phi(z)(y_1,\,0,...,\,0)=y_1^{p^{\a}d}u$$
$$\text{ such that } p\wedge d=1 \text{ and } u \text{ is a unit}\}.$$
By assumption $\Delta$ is not empty. Let us denote by $\a$ the least integer of $\Delta$. Let us choose an element $z_1\in \k\llceil x_1,\,x_2\rrceil[\k']$ such that  $\phi(z_1)=y_1^{p^{\a}d}u$ with $p\wedge d=1$ and  $u$ a unit.
 By using the following modification of type (2): $\psi(y_1)=y_1$ and $\psi(y_k)=y_ky_1$ for $k>1$, we can replace $\phi$ by $\ovl{\phi}:=\psi^k\circ\phi$, for $k\geq 0$, such that  $\ovl{\phi}(z_1)=y_1^{p^{\a}d}u$ for some unit $u$, with $p\wedge d=1$. Let us choose $z_2\in \k\llceil x_1,\,x_2\rrceil[\k']$ such that $(z_1,\,z_2)$ is a regular system of parameters. Now, because $\a$ is the least integer of $\Delta$,  $\ini(\ovl{\phi}(z))$ has no monomial of the form $cy_1^k$ such that $p^{\a}$ does not divide $k$.
 Then we can skip Step 4, and using Step 5,  we can replace $\ovl{\phi}$ by $\ovl{\ovl{\phi}}$ such that $\ovl{\ovl{\phi}}(z_1)=y_1^{p^{\a}de_{1,\,1}}y_2^{p^{\a}de_{1,\,2}}u_1$, $\ovl{\ovl{\phi}}(z_2)= y_1^{i_1e_{1,\,1}+i_2e_{2,\,1}}y_2^{i_1e_{1,\,2}+i_2e_{2,\,2}}u_2$ for some units $u_1$ and $u_2$, and with $\ovl{\phi}=\psi'\circ\ovl{\ovl{\phi}}$ where $\psi'$ is a composition of blowing-ups and automorphisms of $\k\llceil y_1,...,\,y_m\rrceil[\k']$. Moreover the matrix $(e_{i,\,j})_{i,\,j}$ is invertible. 
Then using modifications of type $(4)$ we transform $\ovl{\ovl{\phi}}$ in $\ovl{\ovl{\ovl{\phi}}}$ such that $\ovl{\ovl{\ovl{\phi}}}(x)=y_1^{p^{\a_1}}u_1$ and $\ovl{\ovl{\ovl{\phi}}}(z')=y_2^{p^{\a_2}}u_2$ for some units $u_1$ and $u_2$. Hence $\grk\phi=2$.\\
Now, if $\cha(\k)=0$ then we can do almost the same, but we do not need $\Delta$. We just choose $z_1=x_1$ and $z_2=x_2$. After that the proof is the same as above.
\end{proof}

\begin{corollary}
Let $\phi\ :\ \k\llceil x_1,\,x_2\rrceil\lgw \k\llceil y_1,...,\,y_m\rrceil$ be an injective homomorphism of Henselian $\k$-algebras where $\k$ is infinite. Let $\k\lgw \k'$ be a finite field extension and let us assume that there is a W-system $\k'\llceil x\rrceil$ over $\k'$ such that $\k'\llceil x\rrceil\cap\k[[x]]=\k\llceil x\rrceil$. Then the induced homomorphism  $\phi_{\k'}\ : \ \k'\llceil x_1,\,x_2\rrceil\lgw \k'\llceil y_1,...,\,y_m\rrceil$  is injective.
\end{corollary}

\begin{proof}
By Theorem \ref{dim2} $\grk(\phi)=2$. Thus, if $\cha(\k)=p>0$,  $\phi$ can be transformed using modifications into a homomorphism $\ovl{\phi}$ such that $\ovl{\phi}(x_1)=y_1^{p^{\a_1}}u_1$ and $\ovl{\phi}(x_2)=y_2^{p^{\a_2}}u_2$ for some units. Then $\phi_{\k'}$ can be transformed in the same way and $\grk(\phi_{\k'})=2$. Then $\phi_{\k'}$ is injective by Lemma \ref{lemma4}. The proof in characteristic zero is the same.
\end{proof}

Using this result we deduce the following two results, the first being a generalization to the case of Henselian $\k$-algebras of a theorem of S. S. Abhyankar and M. van der Put (cf. Theorem 2.10 of \cite{A-vdP}):

\begin{theorem}\label{A-vdP}
Let $\phi : A\lgw B$ be a homomorphism of good Henselian  $\k$-algebras where $A$ and $B$ are regular  and $\dim(A)=2$.  If $\phi$ is injective then it is strongly injective.
\end{theorem}

\begin{proof}
We have $\grk(\phi)=2$ from Theorem \ref{dim2}. Hence from Theorem \ref{grk-strongly} $\phi$ is strongly injective by Lemma.
\end{proof}

\begin{corollary}
Let $\phi : A\lgw B$ denote a homomorphism of complete local $\k$-algebras where $A$ is a two dimensional integral domain and $B$ is regular. Then $\phi$ is injective if and only if $\phi$ has a linear Chevalley estimate. 
\end{corollary}

\begin{proof}
It is obvious that $\phi$ is injective if it has a linear Chevalley estimate.\\
 On the other hand the result follows from Theorems \ref{dim2} and \ref{maintheorem}.
\end{proof}

\subsection{The algebraic case}
Here we give a generalization of the main theorem of \cite{To1}, \cite{Be} and \cite{Mi}. The result is the following: \textit{any homomorphism of analytic $\k$-algebras defined by algebraic power series has maximal geometric rank}. This result has been proven for homomorphisms of analytic $\C$-algebras defined by polynomials in the three papers cited above. 
\begin{definition}
Let $\phi\ :\ A\lgw B$ be a homomorphism of local $\k$-algebras.  We define $r_2:=\dim\left(\frac{\wdh{A}}{\Ker(\wdh{\phi})}\right)$ and  $r_3:=\dim\left(\frac{A}{\Ker(\phi)}\right)$. Moreover $r_1:=\grk(\phi)$.
\end{definition}
It is clear that $r_2(\phi)\leq r_3(\phi)$. Moreover, from the definition, we see that $r_1(\wdh{\phi})$ is equal to the geometric rank of the homomorphism $\wdh{A}/\Ker(\wdh{\phi})\lgw \wdh{B}$ induced by $\wdh{\phi}$, and using the Abhyankar's Inequality \cite{Ab} and Lemma \ref{lemma3}  we see that $r_1(\phi)\leq r_2(\phi)$. Thus we always have $r_1(\phi)\leq r_2(\phi)\leq r_3(\phi)$. If $r_1(\phi)=r_2(\phi)$ we say that $\phi$ is \textit{regular}. A difficult theorem of A. Gabrielov asserts that if $\phi : A\lgw B$ is a regular homomorphism when $A$ and $B$ are quotients of convergent power series rings over $\C$ then $r_2(\phi)=r_3(\phi)$, i.e. $\Ker(\wdh{\phi})=\Ker(\phi)\wdh{A}$ (cf. \cite{Ga2}).
\begin{definition} A homomorphism $A\lgw B$ of Henselian $\k$-algebras is said to be a \emph{homomorphism of algebraic $\k$-algebras} if $A$ and $B$ are local $\k$-algebras with respect to the W-system of algebraic power series (Example \ref{ex_W} ii)). 
\end{definition}
\begin{theorem}
Let $\phi\ :\ A\lgw B$ a homomorphism of algebraic $\k$-algebras where $B$ is regular if $\cha(\k)=p>0$. Then $r_1(\phi)=r_3(\phi)$.
\end{theorem}
\begin{proof}
If $\cha(\k)=0$ and $B$ is not regular, by the existence of a resolution of singularities for $\text{Spec}(B)$, there exists a homomorphism of Henselian $\k$-algebras which is a composition of local blow-ups $\psi\ :\ B\lgw \k\langle y_1,...,\,y_m\rangle$. In particular $r_1(\psi\circ\phi)=r_1(\phi)$ by Proposition \ref{compo}. Thus we may assume that $B$ is regular.\\
Let us denote $A':=A/\Ker(\phi)$. Then $d:=\dim(A')=r_3(\phi)$. There exists a finite injective homomorphism $\k\langle x_1,...,\,x_d\rangle\lgw A'$ from the Weierstrass preparation Theorem. Let us denote by $\t$ the homomorphism induced by $\phi$ on $\k\langle x\rangle$. By Lemma \ref{lemma5}, $r_1(\t)=r_1(\phi)$, and because $\t$ is injective, $r_3(\t)=d=r_3(\phi)$.\\
Let $t$ be a variable over $\k$. We may replace $\k$ by $\k(t)$ since Corollary \ref{field_tr} and since the homomorphism induced by $\t$ on $\k(t)\langle x\rangle$ is clearly injective.
Now we appy Theorem \ref{structuretheorem} to $\t$. We get the following commutative diagram:
$$\xymatrix{\k\langle x\rangle \ar[r]^{\t} \ar[d]^{\s_1}& \k\langle y\rangle \ar[d]^{\s_2}\\
\k\langle x\rangle\ar[r]^{\ovl{\t}} & \k\langle y\rangle}$$
where $\ovl{\t}$ is as defined in iii) of Theorem \ref{structuretheorem}. In particular we see that $r_3(\ovl{\t})=r_1(\ovl{\t})$ because $\Ker(\ovl{\t})=(x_{r_1(\t)+1},...,\,x_d)$.
We have $r_1(\t)=r_1(\s_2\circ\t)$ and $r_3(\t)=r_3(\s_2\circ\t)$. Moreover $r_1(\t)=r_1(\ovl{\t})$ according to Corollary \ref{compo}. Thus we only have to prove that $r_3(\t)=r_3(\ovl{\t})$.\\
Let us consider the following commutative diagram:
$$\xymatrix{\k\langle x\rangle \ar[r]^{\t} \ar[d]^{\s}& \k\langle y\rangle \\
\k\langle x\rangle\ar[ru]^{\psi} & }$$
where $\t$ is injective and $\s$ is one of the homomorphisms defined in ii) of Theorem \ref{structuretheorem}. We will prove that $\psi$ is still injective. Thus this will prove by induction that $\ovl{\t}$ is injective and that $r_3(\t)=r_3(\ovl{\t})$.\\
In order to prove that $\psi$ is injective, we have to check the three following cases:
If $\s$ is an isomorphism, then it is clear that $\psi$ is injective.\\
If $\s=\chi_d$ ($d\in\N^*$) is defined by $\chi_d(x_1)=x_1^d$, and $\chi_d(x_i)=x_i$ $\forall i\neq 1$, we can write $d=p^re$ with $e\wedge p=1$. If $f\in \Ker(\psi)$, then let us define $g:=\prod_{\e\in\mathbb{U}_e}\left(f(\e x_1,\,x_2,...,\,x_d)\right)^{p^r}$, where $\mathbb{U}_e$ is the set of $e$-roots of unity in a finite field extension of $\k$. Then $g\in\k\langle x\rangle$ and $g\in\Im(\s)$. Let $g'\in\k\langle x\rangle $ such that $\s(g')=g$. Then $\t(g')=\psi(g)=0$. Thus $g'=0$ because $\t$ is injective, hence $f=0$ and $\psi$ is injective.\\
Finally, let us assume that $\s=q$  defined by $q(x_i)=x_i$ for $i\neq 2$ and $q(x_2)=x_1x_2$. Let $f\in\Ker(\psi)$. Let $P(Y)\in\k[x][Y]$  be an irreducible polynomial having $f$ as a root. Let us denote by $a_i\in\k[x]$, $0\leq i\leq r$, its coefficients (i.e. $a_rf^r+\cdots+a_1f+a_0=0$). Then $\psi(a_0)=0$ and $a_0\neq 0$ if $f\neq 0$. If $d:=\deg_{x_2}(a_0)$, then $g:=x_1^da_0\in\Im(\s)$ and $\psi(g)=0$. Let $g'\in\k[x]$ such that $\s(g')=g$. Then $\t(g')=\psi(g)=0$. Thus $g'=0$ because $\t$ is injective, hence $a_0=0$, then $f=0$ and $\psi$ is injective.\\
\end{proof}

\begin{corollary}\label{To}
Let $\phi\ :\ \k\{x\}\lgw \k\{y\}/I\k\{y\}$ be a homomorphism of analytic $\k$-algebras where $k$ is a valued field, $I$ is an ideal of $\k\langle y\rangle$ and such that $\phi(x_i)\in\k\langle y\rangle/I$ for $1\leq i\leq m$. Assume moreover that $\cha(\k)=0$ or $I=(0)$. Then $r_1(\phi)=r_3(\phi)$.
\end{corollary}
\begin{proof}
Let $\psi\ :\ \k\langle x\rangle\lgw \k\langle y\rangle /I$ be the homomorphism of Henselian $\k$-algebra defined by $\psi(x_i):=\phi(x_i)$ for $1\leq i\leq n$. Then we have $r_1(\psi)=r_3(\psi)$ by the preceding theorem. Moreover, by Lemma \ref{lemma3} , we have $r_1(\psi)=r_1(\wdh{\psi})=r_1(\wdh{\phi})=r_1(\wdh{\phi})$ because $\wdh{\psi}=\wdh{\phi}$. Clearly $\Ker(\psi)\k\{ x\}\subset \Ker(\phi)$, thus $r_3(\phi)\leq r_3(\psi)$. Thus $r_1(\phi)\leq r_3(\phi)\leq r_3(\psi)=r_1(\psi)=r_1(\phi)$ and we get the conclusion.
\end{proof}

\end{document}